\documentclass[11pt]{article}
\usepackage{amsfonts,amscd,amssymb}
\newtheorem{definition}{Definition}[section]

\newtheorem{proposition}[definition]{Proposition}
\newtheorem{corollary}[definition]{Corollary}
\newtheorem{remark}[definition]{Remark}

\newtheorem{theorem}[definition]{Theorem}
\newtheorem{example}[definition]{Example}
\newtheorem{examples}[definition]{Examples}

\def\va{\varepsilon}
\def\v{\varphi}

\def\rh{\rightharpoonup}
\def\lh{\leftharpoonup}

\def\ra{\rightarrow}
\def\a{\alpha}
\def\b{\beta}

\def\l{\lambda}
\def\r{\rho}
\def\cd{\cdot}

\def\O{\Omega}

\def\ov{\overline}

\def\mf{\mathfrak}
\def\mb{\mathbb}
\newcommand{\nat}{\mbox{$\;\natural \;$}}

\newcommand{\mfa}{\mbox{$\mf {a}$}}
\newcommand{\mfb}{\mbox{$\mf {b}$}}

\newcommand{\smi}{\mbox{$S^{-1}$}}
\newcommand{\gsm}{\mbox{$\blacktriangleright \hspace{-0.7mm}<$}}
\newcommand{\gtl}{\mbox{${\;}$$>\hspace{-0.85mm}\blacktriangleleft$${\;}$}}
\newcommand{\trl}{\mbox{${\;}$$\triangleright \hspace{-1.6mm}<$${\;}$}}
\newcommand{\gsl}{\mbox{${\;}$$>\hspace{-1.7mm}\triangleleft$${\;}$}}

\newcommand{\tx}{\mbox{$\tilde {x}$}}
\newcommand{\tX}{\mbox{$\tilde {X}$}}

\newcommand{\tqra}{\mbox{$\tilde {q}^1_{\rho }$}}
\newcommand{\tqrb}{\mbox{$\tilde {q}^2_{\rho }$}}
\newcommand{\tQra}{\mbox{$\tilde {Q}^1_{\rho }$}}
\newcommand{\tQrb}{\mbox{$\tilde {Q}^2_{\rho }$}}
\renewcommand{\theequation}{\thesection.\arabic{equation}}

\def\rawo\lonra{\longrightarrow}

\def\ot{\otimes}

\newcommand{\selabel}[1]{\label{se:#1}}

\newcommand{\eqref}[1]{(\ref{eq:#1})}

\newenvironment{proof}{{\it Proof.}}{\hfill $ \square $ \vskip 4mm}

\begin{document}
\title{L-R-smash product for (quasi) Hopf algebras
\thanks{Research partially supported by the EC programme LIEGRITS, 
RTN 2003, 505078, and by the bilateral project ``New techniques in 
Hopf algebras and graded ring theory'' of the Flemish and Romanian 
Ministries of Research. The first author was also partially supported by  
the programme CERES of the Romanian Ministry of Education and 
Research, contract no. 4-147/2004.}}
\author
{Florin Panaite\\
Institute of Mathematics of the 
Romanian Academy\\ 
PO-Box 1-764, RO-014700 Bucharest, Romania\\
e-mail: Florin.Panaite@imar.ro\\
\and 
Freddy Van Oystaeyen\\
Department of Mathematics and Computer Science\\
University of Antwerp, Middelheimlaan 1\\
B-2020 Antwerp, Belgium\\
e-mail: Francine.Schoeters@ua.ac.be}
\date{}
\maketitle

\begin{abstract}
We introduce a more general version of the so-called L-R-smash product and 
study its relations with other kinds of crossed products (two-sided 
smash and crossed product and diagonal crossed product). We also give an 
interpretation of the L-R-smash product in terms of an   
L-R-twisting datum.  
\end{abstract}
\section*{Introduction}
The L-R-smash product was introduced and studied  
in a series of papers  
\cite{b1}, \cite{b2}, \cite{b3}, \cite{b4}, with motivation and 
examples coming from the theory of deformation quantization. 
It is defined as follows:  
if $H$ is a {\it cocommutative} bialgebra and ${\cal A}$ is an 
$H$-bimodule algebra, 
the L-R-smash product ${\cal A}\nat H$ is an associative  
algebra structure defined on ${\cal A}\ot H$ by the 
multiplication rule 
\begin{eqnarray*}
&&(\varphi \nat h)(\varphi '\nat h')=(\varphi \cdot h'_1)(h_1\cdot \varphi ')
\nat h_2h'_2, \;\;\;\forall \;\varphi , \varphi '\in {\cal A}, \;
h, h'\in H.
\end{eqnarray*}
If the right $H$-action is trivial, ${\cal A}\nat H$ coincides  
with the ordinary smash product ${\cal A}\# H$. \\
Our first remark is that, if we replace the above multiplication by 
\begin{eqnarray*}
&&(\varphi \nat h)(\varphi '\nat h')=(\varphi \cdot h'_2)(h_1\cdot \varphi ')
\nat h_2h'_1, \;\;\;\forall \;\varphi , \varphi '\in {\cal A}, \;
h, h'\in H, 
\end{eqnarray*}
then this multiplication is associative {\it without} the assumption that 
$H$ is cocommutative; this more general object will also be  
denoted by ${\cal A}\nat H$ and called L-R-smash product.\\
Actually, we introduce a much more general construction: we  
define an L-R-smash product ${\cal A}\nat {\mb A}$ (an associative algebra), 
where ${\cal A}$ is a bimodule algebra and ${\mb A}$ is a bicomodule 
algebra over a quasi-bialgebra $H$. Our motivation for defining the 
L-R-smash product in this generality stems from its relations with some  
constructions from \cite{hn1} and \cite{bpvo2}. Indeed, let $H$, ${\cal A}$,  
${\mb A}$ be as above and also $A$ a left $H$-module algebra, $B$ a 
right $H$-module algebra, $\mathfrak{A}$ a right $H$-comodule algebra and 
$\mathfrak{B}$ a left $H$-comodule algebra; then $A\ot B$ becomes an 
$H$-bimodule algebra and $\mf {A}\ot \mf {B}$ an $H$-bicomodule algebra, 
and we prove that we have algebra isomorphisms $(A\ot B)\nat {\mb A}\simeq 
A\gsm {\mb A}\gtl B$ and ${\cal A}\nat (\mf {A}\ot \mf {B})\simeq 
\mf {A}\gsl {\cal A}\trl \mf {B}$, where $A\gsm {\mb A}\gtl B$ and 
$\mf {A}\gsl {\cal A}\trl \mf {B}$ are the two-sided smash and crossed  
products introduced in \cite{bpvo2}, \cite{hn1}. These results combined 
with the ones in \cite{bpvo2} suggest, and we are able to prove, that if 
$H$ is moreover a quasi-Hopf algebra then ${\cal A}\nat {\mb A}$ is 
isomorphic to the generalized diagonal crossed product 
${\cal A}\bowtie {\mb A}$ constructed in \cite{bpvo2} (based on \cite{hn1}). 
As a consequence of this, we get a new realization of the quantum double 
of a finite dimensional quasi-Hopf algebra $H$, having the L-R-smash 
product $H^*\nat H$ as algebra structure. Also as a consequence, we obtain 
that over a cocommutative Hopf algebra $H$, an L-R-smash product 
${\cal A}\nat H$ is actually isomorphic to an ordinary smash product 
${\cal A}\# H$, where the left $H$-action on ${\cal A}$ is now given by 
$h\rightarrow \varphi =h_1\cdot \varphi \cdot S(h_2)$. \\
Note that the multiplication of a diagonal crossed product (such as the 
quantum double of a Hopf or quasi-Hopf algebra) involves the antipode, 
while the one of an L-R-smash product does not. Hence, the L-R-smash 
product is not only a generalization of the ordinary smash product, but 
it can be regarded also as a substitute of the diagonal crossed product for 
bialgebras without antipode. \\
More can be said about the L-R-smash product over bialgebras and Hopf 
algebras. For instance, we provide a Maschke-type theorem for L-R-smash 
products, we give an interpretation (inspired by \cite{fst}) 
of the L-R-smash product in terms of what we call an L-R-twisting datum, 
and we find the counterpart, for twisted products, of the 
isomorphism between the L-R-smash product and the generalized 
diagonal crossed product.\\  
For completeness and possible further use, we also introduce 
the dual version  
of the L-R-smash product (over bialgebras), called L-R-smash coproduct 
(generalizing Molnar's smash coproduct). 
\section{Preliminaries}\selabel{1}
In this section we recall some definitions and results and fix  
notation that will be used throughout the paper.\\[2mm]
We work over a commutative field $k$. All algebras, linear spaces
etc. will be over $k$; unadorned $\ot $ means $\ot_k$. Following
Drinfeld \cite{d1}, a quasi-bialgebra is a fourtuple $(H, \Delta ,
\va , \Phi )$, where $H$ is an associative algebra with unit, 
$\Phi$ is an invertible element in $H\ot H\ot H$, and $\Delta :\
H\ra H\ot H$ and $\va :\ H\ra k$ are algebra homomorphisms
satisfying the identities
\begin{eqnarray}
&&(id \ot \Delta )(\Delta (h))=%
\Phi (\Delta \ot id)(\Delta (h))\Phi ^{-1},\label{q1}\\[1mm]%
&&(id \ot \va )(\Delta (h))=h\ot 1, %
\mbox{${\;\;\;}$}%
(\va \ot id)(\Delta (h))=1\ot h,\label{q2}
\end{eqnarray}
for all $h\in H$, and $\Phi$ has to be a normalized $3$-cocycle,
in the sense that
\begin{eqnarray}
&&(1\ot \Phi)(id\ot \Delta \ot id) (\Phi)(\Phi \ot 1)= (id\ot id
\ot \Delta )(\Phi ) (\Delta \ot id \ot id)(\Phi
),\label{q3}\\[1mm]%
&&(id \ot \va \ot id )(\Phi )=1\ot 1\ot 1.\label{q4}
\end{eqnarray}
The identities (\ref{q2}), (\ref{q3}) and (\ref{q4}) also imply
that
\begin{equation}\label{q7}
(\va \ot id\ot id)(\Phi )= (id \ot id\ot \va )(\Phi )=1\ot 1\ot 1.
\end{equation} 
The map $\Delta $ is called the coproduct or the
comultiplication, $\va $ the counit and $\Phi $ the reassociator.
We use the version of Sweedler's sigma notation: $\Delta (h)=h_1\ot  
h_2$, but since $\Delta$ is only quasi-coassociative we adopt the
further convention%
$$%
(\Delta \ot id)(\Delta (h))= h_{(1, 1)}\ot h_{(1, 2)}\ot h_2, 
\mbox{${\;\;\;}$} (id\ot \Delta )(\Delta (h))=
h_1\ot h_{(2, 1)}\ot h_{(2,2)}, %
$$%
for all $h\in H$. We will denote the tensor components of $\Phi$
by capital letters, and those of $\Phi^{-1}$ by small letters, 
namely
\begin{eqnarray*}
&&\Phi=X^1\ot X^2\ot X^3=T^1\ot T^2\ot T^3=Y^1\ot  
Y^2\ot Y^3=\cdots\\%
&&\Phi^{-1}=x^1\ot x^2\ot x^3=
t^1\ot t^2\ot t^3=y^1\ot y^2\ot y^3=\cdots
\end{eqnarray*}
The quasi-bialgebra $H$ is called a quasi-Hopf algebra if there exists an 
anti-automorphism $S$ of the algebra $H$ and elements $\a , \b \in
H$ such that, for all $h\in H$, we have:
\begin{eqnarray}
&&S(h_1)\a h_2=\va (h)\a \mbox{${\;\;\;}$ and ${\;\;\;}$}
h_1\b S(h_2)=\va (h)\b ,\label{q5}\\[1mm]%
&&X^1\b S(X^2)\a X^3=1 %
\mbox{${\;\;\;}$ and${\;\;\;}$}%
S(x^1)\a x^2\b S(x^3)=1.\label{q6}
\end{eqnarray}
The axioms for a quasi-Hopf algebra imply that $\va (\a )\va (\b 
)=1$, so, by rescaling $\a $ and $\b $, we may assume without loss
of generality that $\va (\a )=\va (\b )=1$ and $\va \circ S=\va $.\\
Next we recall that the definition of a quasi-bialgebra or 
quasi-Hopf algebra is 
"twist covariant" in the following sense. An invertible element
$F\in H\ot H$ is called a {\sl gauge transformation} or {\sl
twist} if $(\va \ot id)(F)=(id\ot \va)(F)=1$. If $H$ is a quasi-bialgebra 
or a quasi-Hopf algebra and $F=F^1\ot F^2\in H\ot H$ is a gauge  
transformation with inverse $F^{-1}=G^1\ot G^2$, then we can 
define a new quasi-bialgebra (respectively quasi-Hopf algebra) 
$H_F$ by keeping the 
multiplication, unit, counit (and antipode in the case of a quasi-Hopf 
algebra) of $H$ and replacing the 
comultiplication, reassociator and the elements $\alpha$ and $\beta$ by 
\begin{eqnarray}
&&\Delta _F(h)=F\Delta (h)F^{-1},\label{g1}\\[1mm]%
&&\Phi_F=(1\ot F)(id \ot \Delta )(F) \Phi (\Delta \ot id)
(F^{-1})(F^{-1}\ot 1),\label{g2}\\[1mm]%
&&\a_F=S(G^1)\a G^2,%
\mbox{${\;\;\;}$}%
\b_F=F^1\b S(F^2).\label{g3}
\end{eqnarray}
It is known that the antipode of a Hopf 
algebra is an anti-coalgebra morphism. For a quasi-Hopf algebra,
we have the following: there exists a gauge 
transformation $f\in H\ot H$ such that
\begin{equation} \label{ca}
f\Delta (S(h))f^{-1}=(S\ot S)(\Delta ^{cop}(h)) 
\mbox{,${\;\;\;}$for all $h\in H$.}
\end{equation}
The element $f$ may be computed explicitly. First set 
\begin{equation}
A^1\ot A^2\ot A^3\ot A^4=(\Phi \ot 1) (\Delta \ot id\ot id)(\Phi ^{-1}),
\end{equation}
\begin{equation} B^1\ot B^2\ot B^3\ot B^4=
(\Delta \ot id\ot id)(\Phi )(\Phi ^{-1}\ot 1),
\end{equation}
and then define $\gamma, \delta\in H\ot H$ by
\begin{equation} \label{gd}%
\gamma =S(A^2)\a A^3\ot S(A^1)\a A^4~~{\rm and}~~ \delta
=B^1\b S(B^4)\ot B^2\b S(B^3).
\end{equation}
Then $f$ and $f^{-1}$ are given by the formulae
\begin{eqnarray}
f&=&(S\ot S)(\Delta ^{cop}(x^1)) \gamma \Delta (x^2\b
S(x^3)),\label{f}\\%
f^{-1}&=&\Delta (S(x^1)\a x^2) \delta (S\ot S)(\Delta
^{cop}(x^3)).\label{g}
\end{eqnarray}
Suppose that $(H, \Delta , \varepsilon , \Phi )$ is a
quasi-bialgebra. If $U,V,W$ are left (right) $H$-modules, define
$a_{U,V,W}, {\bf a}_{U, V, W} :(U\otimes V)\otimes W\rightarrow
U\otimes (V\otimes W)$
by %
\begin{eqnarray*}
&&a_{U,V,W}((u\otimes v)\otimes w)=\Phi \cdot (u\otimes
(v\otimes w)),\\[1mm] %
&&{\bf a}_{U, V, W}((u\ot v)\ot w)= (u\ot (v\ot w))\cd \Phi ^{-1}.
\end{eqnarray*}
The category $_H{\cal M}$ (${\cal M}_H$) of 
left (right) $H$-modules becomes a monoidal category (cf. 
\cite{k}, \cite{m2} for the terminology) with tensor product
$\otimes $ given via $\Delta $, associativity constraints
$a_{U,V,W}$ (${\bf a}_{U, V, W}$), unit $k$ as a trivial
$H$-module and the usual left and right
unit constraints.\\%
Now, let $H$ be a quasi-bialgebra. We say that a $k$-vector space
$A$ is a left $H$-module algebra if it is an algebra in the
monoidal category $_H{\cal M}$, that is $A$ has a multiplication
and a usual unit $1_A$ satisfying the 
following conditions: %
\begin{eqnarray}
&&(a a^{'})a^{''}=(X^1\cd a)[(X^2\cd a^{'})(X^3\cd
a^{''})],\label{ma1}\\[1mm]%
&&h\cd (a a^{'})=(h_1\cd a)(h_2\cd a^{'}),
\label{ma2}\\[1mm]%
&&h\cd 1_A=\va (h)1_A,\label{ma3}
\end{eqnarray}
for all $a, a^{'}, a^{''}\in A$ and $h\in H$, where $h\ot a\ra
h\cd a$ is the left $H$-module structure of $A$. Following
\cite{bpv} we define the smash product $A\# H$ as follows: as
vector space $A\# H$ is $A\ot H$ (elements $a\ot h$ will be
written $a\# h$) with multiplication
given by %
\begin{equation}\label{sm1}
(a\# h)(a^{'}\# h^{'})=%
(x^1\cd a)(x^2h_1\cd a^{'})\# x^3h_2h^{'}, %
\end{equation}
for all $a, a^{'}\in A$, $h, h^{'}\in H$. This $A\# H$ is an
associative algebra with unit $1_A\# 1_H$  
and it is defined by a universal property (as 
Heyneman and Sweedler did for Hopf algebras), see \cite{bpv}.\\
For further use we need also the notion of right $H$-module
algebra. Let $H$ be a quasi-bialgebra. We say that a $k$-linear
space $B$ is a right $H$-module algebra if $B$ is an algebra in
the monoidal category ${\cal M}_H$, i.e. $B$ has a multiplication
and a usual unit $1_B$ satisfying the following conditions:
\begin{eqnarray}
&&(b b^{'})b^{''}=(b\cd x^1)[(b^{'}\cd x^2)(b^{''}\cd x^3)],
\label{rma1}\\[1mm]%
&&(b b^{'})\cd h=(b\cd h_1)(b^{'}\cd h_2),
\label{rma2}\\[1mm]%
&&1_B\cd h=\va (h)1_B,\label{rma3}
\end{eqnarray}
for all $b, b^{'}, b^{''}\in B$ and $h\in H$, where $b\ot h\ra
b\cd h$ is the right $H$-module structure of $B$.\\[2mm]
Recall from \cite{hn1} the notion of comodule algebra over a
quasi-bialgebra.
\begin{definition}
Let $H$ be a quasi-bialgebra. A unital associative algebra
$\mathfrak{A}$ is called a right $H$-comodule algebra if there
exist an algebra morphism $\r :\mathfrak{A}\ra \mathfrak{A}\ot H$
and an invertible element $\Phi _{\r }\in \mathfrak{A}\ot H\ot H$
such that:
\renewcommand{\theequation}{\thesection.\arabic{equation}}
\begin{eqnarray}
&&\Phi _{\r }(\r \ot id)(\r (\mf {a}))=(id\ot \Delta
)(\r (\mf {a}))\Phi _{\r }, 
\mbox{${\;\;\;}$$\forall $ $\mf {a}\in
\mathfrak{A}$,}\label{rca1}\\[1mm]%
&&(1_{\mf {A}}\ot \Phi)(id\ot \Delta \ot id)(\Phi _{\r })(\Phi
_{\r }\ot 1_H)= (id\ot id\ot \Delta )(\Phi _{\r })(\r \ot id\ot
id)(\Phi _{\r }),\label{rca2}\\[1mm]%
&&(id\ot \va)\circ \r =id ,\label{rca3}\\[1mm]%
&&(id\ot \va \ot id)(\Phi _{\r })=(id\ot id\ot \va )(\Phi _{\r }
)=1_{\mathfrak{A}}\ot 1_H.\label{rca4}
\end{eqnarray}
Similarly, a unital associative algebra $\mathfrak{B}$ is called
a left $H$-comodule algebra if there exist an algebra morphism $\l
: \mf {B}\ra H\ot \mathfrak{B}$ and an invertible element $\Phi
_{\l }\in H\ot H\ot \mathfrak{B}$ such that the following
relations hold:
\begin{eqnarray}
&&(id\ot \l )(\l (\mf {b}))\Phi _{\l }=\Phi _{\l
}(\Delta \ot id)(\l (\mf {b})),
\mbox{${\;\;\;}$$\forall $ $\mf {b}\in \mathfrak{B}$,}
\label{lca1}\\[1mm]%
&&(1_H\ot \Phi _{\l })(id\ot \Delta \ot id)(\Phi _{\l })(\Phi \ot
1_{\mf {B}})= (id\ot id\ot \l )(\Phi _{\l })(\Delta \ot id\ot
id)(\Phi _{\l }),\label{lca2}\\[1mm]%
&&(\va \ot id)\circ \l =id ,\label{lca3}\\[1mm]%
&&(id\ot \va \ot id)(\Phi _{\l })=(\va \ot id\ot id)(\Phi _{\l }
)=1_H\ot 1_{\mathfrak{B}}.\label{lca4}
\end{eqnarray}
\end{definition}
When $H$ is a quasi-bialgebra, particular examples of left and
right $H$-comodule algebras are given by $\mf {A}=\mf {B}=H$ and
$\r =\l =\Delta $,
$\Phi _{\r }=\Phi _{\l }=\Phi $.\\ %
For a right $H$-comodule algebra $({\mf A}, \r , \Phi _{\r })$ we
will denote
$$%
\r (\mfa )=\mfa _{<0>}\ot \mfa _{<1>}, \mbox{${\;\;}$} (\r 
\ot id)(\r (\mfa ))=\mfa _{<0, 0>}\ot \mfa _{<0, 1>} \ot \mfa 
_{<1>} \mbox{${\;\;}$etc.}
$$%
for any $\mfa \in {\mf A}$. Similarly, for a left $H$-comodule
algebra $({\mf B}, \l , \Phi _{\l })$, if $\mfb \in {\mf B}$ then
we will denote
$$%
\l (\mfb )=\mfb _{[-1]}\ot \mfb _{[0]}, \mbox{${\;\;}$} 
(id\ot \l )(\l (\mfb ))=\mfb _{[-1]}\ot \mfb _{[0,-1]}\ot 
\mfb _{[0, 0]} \mbox{${\;\;}$etc.}
$$%
In analogy with the notation for the reassociator $\Phi 
$ of $H$, we will write 
\begin{eqnarray*}
&&\Phi _{\r }=\tilde {X}^1_{\r }\ot \tilde {X}^2_{\r }\ot  
\tilde {X}^3_{\r }=
\tilde {Y}^1_{\r }\ot \tilde {Y}^2_{\r }\ot \tilde {Y}^3_{\r }=\cdots  \\
&&\Phi _{\r }^{-1}=\tilde {x}^1_{\r }\ot \tilde {x}^2_{\r }\ot  
\tilde {x}^3_{\r }=\tilde {y}^1_{\r }\ot \tilde {y}^2_{\r 
}\ot \tilde {y}^3_{\r }=\cdots  
\end{eqnarray*}
and similarly for the element $\Phi _{\l }$ of a left $H$-comodule 
algebra $\mf {B}$. \\
If $\mf {A}$ is a right $H$-comodule algebra then we define the elements 
$\tilde{p}_{\r }, \tilde{q}_{\r }\in {\mf A}\ot H$ as follows:
\begin{equation}\label{tpqr}
\tilde {p}_{\r }=\tilde {p}^1_{\r }\ot \tilde {p}^2_{\r}
=\tx ^1_{\r }\ot \tx ^2_{\r }\b S(\tx ^3_{\r }), 
\mbox{${\;\;\;}$}
\tilde {q}_{\r }=\tilde {q}^1_{\r }\ot \tilde {q}^2_{\r}
=\tX ^1_{\r }\ot \smi (\a \tX ^3_{\r })\tX ^2_{\r }.
\end{equation}
By \cite[Lemma 9.1]{hn1}, we have the following relations, for all $\mfa \in
{\mf A}$:
\begin{eqnarray}
&&\hspace*{-1cm}
\r (\mfa_{<0>})\tilde {p}_{\r }[1_{\mf A}\ot S(\mfa _{<1>})]
=\tilde {p}_{\r }[\mfa \ot 1_H],\label{tpqr1}\\%
&&\hspace*{-1cm}
[1_{\mf A}\ot \smi (\mfa_{<1>})]\tilde {q}_{\r }\r
(\mfa_{<0>})=[\mfa \ot 1_H]\tilde {q}_{\r },\label{tpqr1a}\\%
&&\hspace{-1cm}
\r (\tilde {q}^1_{\r })\tilde {p}_{\r }[1_{\mf A}\ot 
S(\tilde{q}^2_{\r })]=1_{\mf A}\ot 1_H\label{tpqr2},\\%
&&\hspace*{-1cm}
[1_{\mf A}\ot \smi (\tilde {p}^2_{\r })]\tilde {q}_{\r }\r
(\tilde {p}^1_{\r })=1_{\mf A}\ot 1_H.\label{tpqr2a} 
\end{eqnarray}
Let $H$ be a quasi-bialgebra, $A$ a left $H$-module algebra and
$\mathfrak{B}$ a left $H$-comodule algebra. Denote by 
$A\gsm \mathfrak{B}$ the $k$-vector space $A\ot\mathfrak{B}$ with newly 
defined multiplication
\begin{equation}\label{gsm}
(a\gsm \mf {b})(a'\gsm \mf {b}')=(\tilde {x}^1_{\l }\cd 
a)(\tilde {x}^2_{\l }\mf {b}_{[-1]}\cd a')\gsm 
\tilde {x}^3_{\l }\mf {b}_{[0]}\mf {b}' ,
\end{equation}
for all $a, a'\in A$ and $\mf {b}, \mf {b}'\in \mathfrak{B}$.
From \cite{bc} it follows that $A\gsm \mathfrak{B}$ is an 
associative algebra 
with unit $1_A\gsm 1_{\mathfrak{B}}$. If we take $\mathfrak{B}=H$ then 
$A\gsm H$ is just the smash product $A\# H$. For this reason the algebra 
$A\gsm \mathfrak{B}$ is called the generalized smash product 
of $A$ and $\mathfrak{B}$. \\[2mm]
The following definition was introduced in \cite{hn1} under the 
name "quasi-commuting pair of $H$-coactions".   
\begin{definition}
Let $H$ be a quasi-bialgebra. By an $H$-bicomodule algebra $\mb {A}$ 
we mean a quintuple $(\l, \r , \Phi _{\l }, \Phi _{\r }, \Phi
_{\l , \r })$, where $\l $ and $\r $ are left and right
$H$-coactions on $\mb {A}$, respectively, and where $\Phi _{\l
}\in H\ot H\ot \mb {A}$, $\Phi _{\r }\in \mb {A}\ot H\ot H$ and
$\Phi _{\l , \r }\in H\ot \mb {A}\ot H$ are invertible elements,
such that:
\begin{itemize}
\item[-]
$(\mb {A}, \l , \Phi _{\l })$ is a left $H$-comodule algebra;
\item[-]$
(\mb {A}, \r , \Phi _{\r })$ is a right $H$-comodule algebra;
\item[-]the following compatibility relations hold:
\renewcommand{\theequation}{\thesection.\arabic{equation}}
\begin{eqnarray}
&&\hspace{-1.8cm}\Phi _{\l , \r }(\l \ot id)(\r (u))=(id\ot \r )(\l 
(u))\Phi _{\l, \r }, \mbox{${\;\;}$$\forall $ $u\in \mb  
{A}$,}\label{bca1}\\[1mm]%
&&\hspace{-1.8cm}(1_H\ot \Phi _{\l , \r })(id\ot \l \ot id)(\Phi
_{\l , \r }) (\Phi _{\l }\ot 1_H)=(id\ot id\ot \r )(\Phi _{\l
})(\Delta \ot
id\ot id)(\Phi _{\l , \r }), \label{bca2}\\[1mm]%
&&\hspace{-1.8cm}(1_H\ot \Phi _{\r })(id\ot \r \ot id)(\Phi _{\l ,\r
})(\Phi _{\l , \r }\ot 1_H)= (id\ot id\ot \Delta )(\Phi _{\l , \r
})(\l \ot id\ot id) (\Phi _{\r }).\label{bca3}
\end{eqnarray}
\end{itemize}
\end{definition}
As pointed out in \cite{hn1}, if $\mb {A}$ is a bicomodule algebra
then, in addition, we have that
\renewcommand{\theequation}{\thesection.\arabic{equation}}
\begin{equation}\label{bca4}
(id_H\ot id_{\mb {A}}\ot \va )(\Phi _{\l , \r })=1_H\ot 1_{\mb
{A}}, \mbox{${\;\;}$} (\va \ot id_{\mb {A}}\ot id_H)(\Phi _{\l ,
\r })= 1_{\mb {A}} \ot 1_H.
\end{equation}
A first example of a bicomodule algebra is $\mb {A}=H$, $\l =\r  
=\Delta $ and $\Phi _{\l }=\Phi _{\r }= \Phi _{\l , \r }=\Phi $.
Related to the left and right comodule algebra structures of $\mb
{A}$ we keep notation as above. For  
simplicity we denote
\begin{eqnarray*}
&&\Phi _{\l , \r }=\Theta ^1\ot \Theta ^2\ot \Theta ^3=
{\bf \Theta }^1\ot {\bf \Theta }^2\ot {\bf \Theta }^3=
\ov {\Theta }^1\ot \ov {\Theta }^2\ot \ov {\Theta }^3, \\
&&\Phi ^{-1}_{\l , \r }=\theta ^1\ot \theta ^2\ot \theta
^3=\tilde{\theta }^1\ot \tilde{\theta }^2\ot \tilde{\theta }^3=
\ov {\theta }^1\ot \ov {\theta }^2\ot \ov {\theta }^3. 
\end{eqnarray*}
Let us denote by $_H{\cal M}_H$ the category of $H$-bimodules; it is 
also a monoidal category, the 
associativity constraints being given by ${\bf a^{'}}_{U, V, W}:
(U\ot V)\ot W\ra U\ot (V\ot W)$,
\renewcommand{\theequation}{\thesection.\arabic{equation}}
\begin{equation}\label{bim}
{\bf a^{'}}_{U, V, W}((u\ot v)\ot w)= \Phi \cd (u\ot (v\ot w))\cd
\Phi ^{-1}, 
\end{equation}
for any $U, V, W\in {}_H{\cal M}_H$ and $u\in U$, $v\in V$ and $w\in
W$. Therefore, we may define algebras in the category of 
$H$-bimodules. Such an algebra will be called an $H$-bimodule
algebra. More exactly, a $k$-vector space ${\cal A}$ is an
$H$-bimodule algebra if ${\cal A}$ is an $H$-bimodule (denote
the actions by $h\cd \varphi $ and $\varphi \cd h$, for $h\in H$ and 
$\varphi \in {\cal A}$)   
which has a multiplication and a usual unit $1_{\cal A}$ such that the
following relations hold:
\renewcommand{\theequation}{\thesection.\arabic{equation}}
\begin{eqnarray}
&&\hspace{-2cm}(\varphi \psi )\xi =(X^1\cd \varphi \cd
x^1)[(X^2\cd \psi \cd x^2)(X^3\cd \xi \cd x^3)],
\mbox{${\;\;}$$\forall $ $\varphi , \psi , \xi \in
{\cal A}$,}\label{bma1}\\[1mm]%
&&\hspace{-2cm}h\cd (\varphi \psi)=(h_1\cd \varphi)(h_2\cd
\psi ), 
\mbox{${\;\;}$} 
(\varphi \psi )\cd h=(\varphi \cd h_1)(\psi \cd h_2),  
\mbox{${\;\;}$$\forall $ $\varphi , \psi \in {\cal A}$, $h\in H$,} 
\label{bma2}\\
&&\hspace{-2cm}
h\cd 1_{\cal A}=\va (h)1_{\cal A},
\mbox{${\;\;}$}
1_{\cal A}\cd h=\va (h)1_{\cal A},
\mbox{${\;\;}$$\forall \; h\in H$.}\label{bma3}
\end{eqnarray}
Let $H$ be a quasi-bialgebra. Then $H^*$, 
the linear dual of $H$, is an $H$-bimodule via the 
$H$-actions %
\renewcommand{\theequation}{\thesection.\arabic{equation}}
\begin{equation}\label{dhbima}
<h\rh \v , h^{'}>=\v (h^{'}h), \mbox{${\;\;\;}$}
<\v \lh h, h^{'}>=\v (hh^{'}),
\mbox{${\;\;\;}$$\forall $ $\varphi \in H^*$, $h, h^{'}\in H$.}
\end{equation}
The convolution $<\v \psi , h>=\sum \v (h_1)\psi (h_2)$, $\v ,
\psi \in H^*$, $h\in H$, is a multiplication on $H^*$; it is not
in general associative, but with this multiplication $H^*$ becomes
an $H$-bimodule algebra.
\section{L-R-smash product over quasi-bialgebras and\\ 
quasi-Hopf algebras}\label{sectiune}
\setcounter{equation}{0}
We introduce the general version of the L-R-smash product as follows. 
\begin{proposition}
Let $H$ be a quasi-bialgebra, ${\cal A}$ an $H$-bimodule 
algebra and ${\mb A}$ an $H$-bicomodule algebra. Define on ${\cal A} 
\ot {\mb A}$ the product  
\begin{eqnarray}\label{nat}
&&(\varphi \nat u)(\psi \nat u^{'})=(\tilde{x}^1_{\lambda }\cd \varphi \cd 
\theta ^3u'_{<1>}\tilde{x}^2_{\rho })
(\tilde{x}^2_{\lambda }u_{[-1]}\theta ^1\cd \psi \cd 
\tilde{x}^3_{\rho })\nat \tilde{x}^3_{\lambda }u_{[0]}\theta ^2u'_{<0>}
\tilde{x}^1_{\rho }
\end{eqnarray}
for $\v , \psi \in {\cal A}$ and $u, u^{'}\in {\mb A}$,  
where $\Phi _{\r }^{-1}=\tx ^1_{\r }\ot \tx ^2_{\r }\ot \tx ^3_{\r  
}$, $\Phi _{\l }^{-1}=\tx ^1_{\l }\ot \tx ^2_{\l }\ot \tx  
^3_{\l }$,  
$\Phi _{\l , \r }^{-1}=\theta ^1\ot \theta ^2\ot  
\theta ^3$,    
and we write $\v \nat u$ in place of $\v \ot u$      
to distinguish the new algebraic structure.    
Then this product defines on 
${\cal A}\ot {\mb A}$ a structure of associative algebra with unit 
$1_{\cal A}\nat 1_{\mb A}$, denoted by     
${\cal A}\nat {\mb A}$ and called the L-R-smash product.
\end{proposition}
\begin{proof}
For $\v , \psi , \xi \in {\cal A}$ and $u, u^{'}, u^{''}\in {\mb A}$ 
we compute:\\[2mm]%
${\;\;\;\;\;\;\;\;\;\;\;}$%
$[(\v \nat u)(\psi \nat u^{'})](\xi \nat u^{''})$%
\begin{eqnarray*}
{\rm (\ref{nat})}&=&[(\tilde{x}^1_{\lambda }\cd \varphi \cd 
\theta ^3u'_{<1>}\tilde{x}^2_{\rho })
(\tilde{x}^2_{\lambda }u_{[-1]}\theta ^1\cd \psi \cd 
\tilde{x}^3_{\rho })\nat \tilde{x}^3_{\lambda }u_{[0]}\theta ^2u'_{<0>}
\tilde{x}^1_{\rho }](\xi \nat u'')\\
{\rm (\ref{nat})}&=&[((\tilde{y}^1_{\lambda })_1\tilde{x}^1_{\lambda }
\cd \varphi \cd 
\theta ^3u'_{<1>}\tilde{x}^2_{\rho }\overline{\theta }^3_1u''_{<1>_1}
(\tilde{y}^2_{\rho })_1)\\
&&((\tilde{y}^1_{\lambda })_2\tilde{x}^2_{\lambda }u_{[-1]}\theta ^1\cd \psi  
\cd \tilde{x}^3_{\rho }\overline{\theta }^3_2u''_{<1>_2}
(\tilde{y}^2_{\rho })_2)]\\
&&(\tilde{y}^2_{\lambda }(\tilde{x}^3_{\lambda })_{[-1]}(u_{[0]})_{[-1]}
\theta ^2_{[-1]}(u'_{<0>})_{[-1]}(\tilde{x}^1_{\rho })_{[-1]}
\overline{\theta }^1\cdot \xi \cdot \tilde{y}^3_{\rho })\\
&&\nat \tilde{y}^3_{\lambda }(\tilde{x}^3_{\lambda })_{[0]}
(u_{[0]})_{[0]}\theta ^2_{[0]}(u'_{<0>})_{[0]}
(\tilde{x}^1_{\rho })_{[0]}\overline{\theta }^2u''_{<0>}\tilde{y}^1_{\rho }\\
{\rm (\ref{lca2}, \ref{bca3})}&=&
[(t^1\tilde{y}^1_{\lambda }
\cd \varphi \cd 
\theta ^3u'_{<1>}\tilde{\theta }^3(\overline{\theta }^2)_{<1>}
\tilde{x}^2_{\rho }u''_{<1>_1}
(\tilde{y}^2_{\rho })_1)\\
&&(t^2(\tilde{y}^2_{\lambda })_1\tilde{x}^1_{\lambda }u_{[-1]}
\theta ^1\cd \psi   
\cd \overline{\theta }^3\tilde{x}^3_{\rho }u''_{<1>_2} 
(\tilde{y}^2_{\rho })_2)]\\
&&(t^3(\tilde{y}^2_{\lambda })_2\tilde{x}^2_{\lambda }(u_{[0]})_{[-1]} 
\theta ^2_{[-1]}(u'_{<0>})_{[-1]}\tilde{\theta }^1
\overline{\theta }^1\cdot \xi \cdot \tilde{y}^3_{\rho })\\ 
&&\nat \tilde{y}^3_{\lambda }\tilde{x}^3_{\lambda }
(u_{[0]})_{[0]}\theta ^2_{[0]}(u'_{<0>})_{[0]}\tilde{\theta }^2
(\overline{\theta }^2)_{<0>}\tilde{x}^1_{\rho }u''_{<0>}\tilde{y}^1_{\rho }\\
{\rm (\ref{bma1}, \ref{rca1}, \ref{lca1})}&=&
(\tilde{y}^1_{\lambda }\cd \varphi \cd  
\theta ^3u'_{<1>}\tilde{\theta }^3(\overline{\theta }^2)_{<1>}
(u''_{<0>})_{<1>}\tilde{x}^2_{\rho }
(\tilde{y}^2_{\rho })_1t^1)\\
&&[((\tilde{y}^2_{\lambda })_1u_{[-1]_1}\tilde{x}^1_{\lambda }
\theta ^1\cd \psi   
\cd \overline{\theta }^3u''_{<1>}\tilde{x}^3_{\rho } 
(\tilde{y}^2_{\rho })_2t^2)\\
&&((\tilde{y}^2_{\lambda })_2u_{[-1]_2}\tilde{x}^2_{\lambda } 
\theta ^2_{[-1]}(u'_{<0>})_{[-1]}\tilde{\theta }^1
\overline{\theta }^1\cdot \xi \cdot \tilde{y}^3_{\rho }t^3)]\\ 
&&\nat \tilde{y}^3_{\lambda }u_{[0]}\tilde{x}^3_{\lambda }
\theta ^2_{[0]}(u'_{<0>})_{[0]}\tilde{\theta }^2
(\overline{\theta }^2)_{<0>}(u''_{<0>})_{<0>}\tilde{x}^1_{\rho }
\tilde{y}^1_{\rho }\\
{\rm (\ref{bca1})}&=&
(\tilde{y}^1_{\lambda }\cd \varphi \cd  
\theta ^3\tilde{\theta }^3(u'_{[0]})_{<1>}(\overline{\theta }^2)_{<1>}
(u''_{<0>})_{<1>}\tilde{x}^2_{\rho }
(\tilde{y}^2_{\rho })_1t^1)\\
&&[((\tilde{y}^2_{\lambda })_1u_{[-1]_1}\tilde{x}^1_{\lambda }
\theta ^1\cd \psi   
\cd \overline{\theta }^3u''_{<1>}\tilde{x}^3_{\rho } 
(\tilde{y}^2_{\rho })_2t^2)\\
&&((\tilde{y}^2_{\lambda })_2u_{[-1]_2}\tilde{x}^2_{\lambda } 
\theta ^2_{[-1]}\tilde{\theta }^1u'_{[-1]}
\overline{\theta }^1\cdot \xi \cdot \tilde{y}^3_{\rho }t^3)]\\ 
&&\nat \tilde{y}^3_{\lambda }u_{[0]}\tilde{x}^3_{\lambda }
\theta ^2_{[0]}\tilde{\theta }^2(u'_{[0]})_{<0>}
(\overline{\theta }^2)_{<0>}(u''_{<0>})_{<0>}\tilde{x}^1_{\rho }
\tilde{y}^1_{\rho }\\
{\rm (\ref{rca2}, \ref{bca2})}&=&
(\tilde{y}^1_{\lambda }\cd \varphi \cd  
\theta ^3(\tilde{x}^3_{\lambda })_{<1>}(u'_{[0]})_{<1>}
(\overline{\theta }^2)_{<1>}
(u''_{<0>})_{<1>}(\tilde{x}^1_{\rho })_{<1>}
\tilde{y}^2_{\rho })\\
&&[((\tilde{y}^2_{\lambda })_1u_{[-1]_1}\theta ^1_1\tilde{x}^1_{\lambda } 
\cd \psi    
\cd \overline{\theta }^3u''_{<1>}\tilde{x}^2_{\rho }  
(\tilde{y}^3_{\rho })_1)\\
&&((\tilde{y}^2_{\lambda })_2u_{[-1]_2}\theta ^1_2\tilde{x}^2_{\lambda } 
u'_{[-1]}
\overline{\theta }^1\cdot \xi \cdot \tilde{x}^3_{\rho }
(\tilde{y}^3_{\rho })_2)]\\ 
&&\nat \tilde{y}^3_{\lambda }u_{[0]}\theta ^2(\tilde{x}^3_{\lambda })_{<0>}
(u'_{[0]})_{<0>}
(\overline{\theta }^2)_{<0>}(u''_{<0>})_{<0>}(\tilde{x}^1_{\rho })_{<0>}
\tilde{y}^1_{\rho }\\
{\rm (\ref{nat})}&=&(\varphi \nat u)[(\tilde{x}^1_{\lambda }\cdot \psi \cdot  
\overline{\theta }^3u''_{<1>}\tilde{x}^2_{\rho })(\tilde{x}^2_{\lambda }
u'_{[-1]}\overline{\theta }^1\cdot \xi \cdot \tilde{x}^3_{\rho })\nat 
\tilde{x}^3_{\lambda }u'_{[0]}\overline{\theta }^2u''_{<0>}
\tilde{x}^1_{\rho }]\\
{\rm (\ref{nat})}&=&(\varphi \nat u)[(\psi \nat u')(\xi \nat u'')],
\end{eqnarray*}
hence the multiplication is associative. It is easy to check that 
$1_{\cal A}\nat 1_{\mb A}$ is the unit.
\end{proof}
\begin{remark} \rm
It is easy to see that, in ${\cal A}\nat {\mb A}$, we have 
$(1\nat u)(1\nat u')=1\nat uu'$ for all $u, u'\in {\mb A}$, hence the map 
${\mb A}\rightarrow {\cal A}\nat {\mb A}$, $u\mapsto 1\nat u$, is an algebra 
map, and     
$(\varphi \nat 1)(1\nat u)=\varphi \cdot u_{<1>}\nat u_{<0>}$. 
\end{remark}  
\begin{examples}\rm 
${\;\;}$1) Let $A$ be a left $H$-module algebra. Then $A$ becomes   
an $H$-bimodule algebra, 
with right $H$-action given via $\va $. In this  
case the multiplication of $A\nat {\mb A}$ becomes 
\begin{eqnarray*} 
&&(a\nat u)(a'\nat u^{'})=(\tilde{x}^1_{\lambda }\cd a)
(\tilde{x}^2_{\lambda }u_{[-1]}\cd \psi )
\nat \tilde{x}^3_{\lambda }u_{[0]}u', 
\end{eqnarray*}
for all $a, a'\in A$ and $u, u'\in {\mb A}$, hence in this case 
$A\nat {\mb A}$ coincides with the generalized smash product 
$A\gsm {\mb A}$.\\
${\;\;}$2) As we have already mentioned, $H$ itself 
is an $H$-bicomodule algebra. So, 
in this case,  the multiplication of 
${\cal A}\nat H$ specializes to     
\begin{eqnarray}
&&(\varphi \nat h)(\psi \nat h^{'})= 
(x^1\cd \varphi \cd 
t^3h'_2y^2)
(x^2h_1t^1\cd \psi \cd  
y^3)\nat x^3h_2t^2h'_1y^1, 
\end{eqnarray}
for all $\varphi , \psi \in {\cal A}$ and $h, h^{'}\in H$. If the right  
$H$-module structure of ${\cal A}$ is trivial, then ${\cal A}\nat H$ 
coincides with the smash product ${\cal A}\# H$.\\ 
${\;\;}$3) Let $H$ be an ordinary bialgebra,  
${\cal A}$ an $H$-bimodule  
algebra and $\mb {A}$ an $H$-bicomodule algebra in the usual (Hopf) sense.  
In this  
case the multiplication of ${\cal A}\nat \mb {A}$ becomes:
\begin{eqnarray}
&&(\varphi \nat u)(\psi \nat u^{'})=  
(\varphi \cd u'_{<1>})(u_{[-1]}\cdot \psi )
\nat u_{[0]}u'_{<0>},
\end{eqnarray} 
for all $\varphi , \psi \in {\cal A}$ and $u, u^{'}\in \mb {A}$. If 
moreover ${\mb A}=H$, the multiplication of ${\cal A}\nat H$ is  
\begin{eqnarray}
&&(\varphi \nat h)(\psi \nat h')=  
(\varphi \cd h'_2)(h_1\cdot \psi )
\nat h_2h'_1,\;\;\;\forall \;\varphi , \psi \in {\cal A}, \;h, h'\in H.
\end{eqnarray}
In case $H$ is cocommutative, this product can be written as 
\begin{eqnarray}
&&(\varphi \nat h)(\psi \nat h')=  
(\varphi \cd h'_1)(h_1\cdot \psi )
\nat h_2h'_2,
\end{eqnarray}
and this was the original L-R-smash product defined in \cite{b1}, \cite{b2}, 
\cite{b3}, \cite{b4}.  
\end{examples}
Recall from \cite{bpvo2} the so-called two-sided generalized smash  
product, defined as follows. 
Let $H$ be a quasi-bialgebra, $A$ a left $H$-module algebra, $B$ a 
right $H$-module algebra and $\mb A$ an $H$-bicomodule algebra. If
we define on $A\ot {\mb A}\ot B$ a multiplication, by\\[2mm]
${\;\;\;}$%
$(a\gsm u\gtl b)(a^{'}\gsm u^{'}\gtl b^{'})$
\begin{equation}\label{tgsm}
=(\tx ^1_{\l }\cd a)(\tx ^2_{\l }u_{[-1]}\theta ^1\cd
a^{'})\gsm \tx ^3_{\l } u_{[0]}\theta ^2u^{'}_{<0>}\tx ^1_{\r
}\gtl (b\cd \theta ^3u^{'}_{<1>}\tx ^2_{\r })(b^{'}\cd \tx ^3_{\r}), 
\end{equation}
for all $a, a^{'}\in A$, $u, u^{'}\in {\mb A}$ and $b, b^{'}\in
B$ (we write $a\gsm u\gtl b$ for $a\ot u\ot b$), and we
denote this structure by $A\gsm {\mb A}\gtl 
B$, then it is an associative algebra with unit 
$1_A\gsm 1_{\mb A} \gtl 1_B$.\\
Note that, given $A, B$ as above, $A\ot B$ becomes an $H$-bimodule algebra, 
with $H$-actions  
\begin{equation}\label{bmas}
h\cd (a\ot b)\cd h^{'}=h\cd a\ot b\cd h^{'}, 
\mbox{${\;\;\;}$$\forall $ $a\in A$, $h, h^{'}\in H$, $b\in B$.}
\end{equation}
\begin{proposition}\label{two-sided}
If $H, A, B, {\mb A}$ are as above, then we have an algebra isomorphism 
\begin{eqnarray*}
&&\phi :(A\ot B)\nat {\mb A}\simeq A\gsm {\mb A}\gtl B, \\
&&\phi ((a\ot b)\nat u)=a\gsm u\gtl b,\;\;\forall \;\;a\in A, \;b\in B,\;
u\in {\mb A}.
\end{eqnarray*}
\end{proposition}
\begin{proof} We compute:\\[2mm]
${\;\;\;\;\;}$
$\phi ([(a\ot b)\nat u)][(a'\ot b')\nat u'])$
\begin{eqnarray*}
&&=\phi ((\tilde{x}^1_{\lambda }\cd (a\ot b)\cd 
\theta ^3u'_{<1>}\tilde{x}^2_{\rho })
(\tilde{x}^2_{\lambda }u_{[-1]}\theta ^1\cd (a'\ot b') \cdot 
\tilde{x}^3_{\rho })\nat \tilde{x}^3_{\lambda }u_{[0]}\theta ^2u'_{<0>}
\tilde{x}^1_{\rho })\\
&&=\phi ((\tilde{x}^1_{\lambda }\cd a\ot b\cd    
\theta ^3u'_{<1>}\tilde{x}^2_{\rho })
(\tilde{x}^2_{\lambda }u_{[-1]}\theta ^1\cd a'\ot b' \cdot  
\tilde{x}^3_{\rho })\nat \tilde{x}^3_{\lambda }u_{[0]}\theta ^2u'_{<0>}
\tilde{x}^1_{\rho })\\
&&=\phi (((\tilde{x}^1_{\lambda }\cdot a)(\tilde{x}^2_{\lambda }u_{[-1]}
\theta ^1\cd a')\ot (b\cd \theta ^3u'_{<1>}\tilde{x}^2_{\rho })
(b'\cdot \tilde{x}^3_{\rho }))\nat \tilde{x}^3_{\lambda }u_{[0]}\theta ^2
u'_{<0>}\tilde{x}^1_{\rho })\\
&&=(\tilde{x}^1_{\lambda }\cdot a)(\tilde{x}^2_{\lambda }u_{[-1]}
\theta ^1\cd a')\gsm \tilde{x}^3_{\lambda }u_{[0]}\theta ^2
u'_{<0>}\tilde{x}^1_{\rho }\gtl (b\cd \theta ^3u'_{<1>}\tilde{x}^2_{\rho })
(b'\cdot \tilde{x}^3_{\rho })\\
&&=(a\gsm u\gtl b)(a'\gsm u'\gtl b')\\
&&=\phi ((a\ot b)\nat u)\phi ((a'\ot b')\nat u'),
\end{eqnarray*} 
and the proof is finished.
\end{proof}
Recall from \cite{bpvo2} the so-called generalized two-sided  
crossed product, defined as follows: if $H$ is  
a quasi-bialgebra, $\mf {A}$ a right $H$-comodule algebra, $\mf {B}$ a 
left $H$-comodule algebra and ${\cal A}$ an $H$-bimodule algebra, define on  
$\mf {A}\ot {\cal A}\ot \mf {B}$ a multiplication by the formula 
\begin{eqnarray}
&&\hspace*{-2cm} 
(\mf {a}\gsl \varphi \trl 
\mf {b})(\mf {a}'\gsl \varphi ' \trl \mf {b}')\nonumber\\
&=&\mf {a}\mf {a}'_{<0>}\tx ^1_{\r }\gsl (\tx 
^1_{\l }\cdot \varphi \cdot   
\mf {a}'_{<1>}\tx ^2_{\r })(\tx ^2_{\l }\mf {b}_{[-1]}\cdot  \varphi ' 
\cdot \tx ^3_{\r })\trl \tx ^3_{\l }\mf {b}_{[0]}\mf {b}' ,\label{gtscp}
\end{eqnarray}
for all $\mf {a}, \mf {a}'\in \mf {A}$, $\mf {b}, \mf
{b}'\in \mf {B}$ and $\v, \v ' \in {\cal A}$,  
where we write $\mf {a}\gsl \varphi \trl 
\mf {b}$ for $\mf {a}\ot \varphi \ot \mf {b}$. Then  
this multiplication yields an associative algebra structure  
with unit $1_{\mf {A}}\gsl 1_{{\cal A}} \trl 1_{\mf {B}}$, denoted  
by $\mf {A}\gsl {\cal A}\trl \mf {B}$. For  
$H$ finite dimensional and ${\cal A}=H^*$ we recover the two-sided 
crossed product $\mf {A}\gsl H^*\trl \mf {B}$ from \cite{hn1}.\\
Note that, given $\mf {A}$, $\mf {B}$ as above,  
$\mf {A}\ot \mf {B}$ becomes an $H$-bicomodule algebra, with the following   
structure:  
$\r (\mfa \ot \mfb)=(\mfa _{<0>}\ot \mfb )\ot \mfa _{<1>}$, 
$\l (\mfa \ot \mfb )=\mfb _{[-1]}\ot (\mfa \ot \mfb _{[0]})$, 
$\Phi _{\r }=(\tilde {X}^1_{\r }\ot 1_{\mf {B}})\ot  
\tilde {X}^2_{\r }\ot \tilde {X}^3_{\r }$, 
$\Phi _{\l }=\tilde {X}^1_{\l }\ot    
\tilde {X}^2_{\l }\ot (1_{\mf {A}}\ot \tilde {X}^3_{\l })$, 
$\Phi _{\l , \r }=1_H\ot (1_{\mf {A}}\ot 1_{\mf {B}})\ot 1_H$, 
for all $\mfa \in \mf {A}$ and $\mfb \in \mf {B}$, see \cite{hn1}. 
\begin{proposition}\label{crossed}
If $H$, ${\cal A}$, $\mf {A}$, $\mf {B}$ are  
as above, then we have an algebra isomorphism  
\begin{eqnarray*}
&&\tau :{\cal A}\nat (\mf {A}\ot \mf {B})\simeq 
\mf {A}\gsl {\cal A}\trl \mf {B}, \\
&&\tau (\varphi \nat (\mathfrak{a}\ot \mathfrak{b}))=\mathfrak{a}\gsl  
\varphi \trl \mathfrak{b}, \;\;\;\forall \;\varphi \in {\cal A},\; 
\mathfrak{a}\in \mf {A}, \;\mathfrak{b}\in \mf{B}. 
\end{eqnarray*}
\end{proposition}
\begin{proof}
We compute:\\[2mm]
${\;\;\;\;}$$\tau ((\varphi \nat (\mathfrak{a}\ot 
\mathfrak{b}))(\varphi '\nat (\mathfrak{a}'\ot 
\mathfrak{b}')))$
\begin{eqnarray*}
&&=\tau ((\tilde{x}^1_{\lambda }\cdot \varphi \cdot (\mathfrak{a}'\ot 
\mathfrak{b}')_{<1>}\tilde{x}^2_{\rho })(\tilde{x}^2_{\lambda }(\mathfrak{a}
\ot \mathfrak{b})_{[-1]}\cdot \varphi '\cdot \tilde{x}^3_{\rho })\\
&&\;\;\;\;\;\;\;\;\nat  
(1_{\mf {A}}\ot \tilde{x}^3_{\lambda })(\mathfrak{a}\ot \mathfrak{b})_{[0]}
(\mathfrak{a}'\ot \mathfrak{b}')_{<0>}(\tilde{x}^1_{\rho }\ot 1_{\mf {B}}))\\
&&=\tau ((\tilde{x}^1_{\lambda }\cdot \varphi \cdot \mathfrak{a}' 
_{<1>}\tilde{x}^2_{\rho })(\tilde{x}^2_{\lambda }
\mathfrak{b}_{[-1]}\cdot \varphi '\cdot \tilde{x}^3_{\rho })\nat 
(1_{\mf {A}}\ot \tilde{x}^3_{\lambda })(\mathfrak{a}\ot \mathfrak{b}_{[0]})
(\mathfrak{a}'_{<0>}\ot \mathfrak{b}')(\tilde{x}^1_{\rho }\ot 1_{\mf {B}}))\\
&&=\tau ((\tilde{x}^1_{\lambda }\cdot \varphi \cdot \mathfrak{a}' 
_{<1>}\tilde{x}^2_{\rho })(\tilde{x}^2_{\lambda }
\mathfrak{b}_{[-1]}\cdot \varphi '\cdot \tilde{x}^3_{\rho })\nat 
(\mathfrak{a}\mathfrak{a}'_{<0>}\tilde{x}^1_{\rho }\ot \tilde{x}^3_{\lambda }
\mathfrak{b}_{[0]}\mathfrak{b}'))\\
&&=\mathfrak{a}\mathfrak{a}'_{<0>}\tilde{x}^1_{\rho }\gsl  
(\tilde{x}^1_{\lambda }\cdot \varphi \cdot \mathfrak{a}' 
_{<1>}\tilde{x}^2_{\rho })(\tilde{x}^2_{\lambda }
\mathfrak{b}_{[-1]}\cdot \varphi '\cdot \tilde{x}^3_{\rho })\trl 
\tilde{x}^3_{\lambda }\mathfrak{b}_{[0]}\mathfrak{b}'\\
&&=(\mathfrak{a}\gsl \varphi \trl \mathfrak{b})
(\mathfrak{a}'\gsl \varphi '\trl \mathfrak{b}')\\
&&=\tau (\varphi \nat (\mathfrak{a}\ot \mathfrak{b}))
\tau (\varphi '\nat (\mathfrak{a}'\ot \mathfrak{b}')),
\end{eqnarray*}
and the proof is finished.
\end{proof}   
Assume that H is a quasi-Hopf algebra, ${\cal A}$ is an $H$-bimodule  
algebra and ${\mb A}$ an $H$-bicomodule algebra. Recall from \cite{bpvo2}  
the so-called generalized diagonal crossed product (which for ${\cal A}=H^*$ 
gives the diagonal crossed product from \cite{hn1}).  
Namely, define the element  
\begin{eqnarray}
&&
\hspace*{-1cm}
\O =(\tX ^1_{\r })_{[-1]_1}\tx ^1_{\l }\theta ^1\ot (\tX
^1_{\r })_{[-1]_2}\tx ^2_{\l }\theta ^2_{[-1]}\ot (\tX ^1_{\r
})_{[0]}\tx ^3_{\l }\theta ^2_{[0]}\ot \smi (f^1\tX ^2_{\r }\theta
^3)\ot \smi (f^2\tX ^3_{\r }) \label{o}
\end{eqnarray}
in $H^{\ot 2}\ot {\mb A}\ot H^{\ot 2}$, 
where $\Phi _{\r }=\tX ^1_{\r }\ot \tX ^2_{\r }\ot \tX ^3_{\r 
}$, $\Phi _{\l }^{-1}=\tx ^1_{\l }\ot \tx ^2_{\l }\ot \tx 
^3_{\l }$,  
$\Phi _{\l , \r }^{-1}=\theta ^1\ot \theta ^2\ot  
\theta ^3$ and $f=f^1\ot f^2$ is the twist defined in (\ref{f}). Then  
define a multiplication on ${\cal A}\ot {\mb A}$, by: 
\begin{eqnarray}\label{gdp}
(\varphi \bowtie u)(\psi \bowtie u^{'})=(\O ^1\cd \varphi \cd  
\O ^5)(\O ^2u_{<0>_{[-1]}}\cd \psi \cd \smi (u_{<1>})\O ^4)\bowtie
\O ^3u_{<0>_{[0]}}u^{'}, 
\end{eqnarray}
for all $\varphi , \psi \in {\cal A}$ and $u, u'\in {\mb A}$. Then this  
multiplication defines an associative algebra structure with unit  
$1_{\cal A}\bowtie 1_{\mb A}$, which will be denoted by 
${\cal A}\bowtie {\mb A}$. \\
It was proved in \cite{bpvo2} that, for an $H$-bimodule algebra of type 
$A\ot B$, where $A$ is a left $H$-module algebra and $B$ is a right 
$H$-module algebra, we have $(A\ot B)\bowtie {\mb A}\simeq A\gsm {\mb A} \gtl 
B$, hence, by Proposition \ref{two-sided}, we obtain 
$(A\ot B)\bowtie {\mb A}\simeq (A\ot B)\nat {\mb A}$. Also, it was proved 
in \cite{bpvo2} that, for an $H$-bicomodule algebra of type 
$\mf {A}\ot \mf {B}$, where $\mf A$ is a right $H$-comodule algebra and 
$\mf {B}$ is a left $H$-comodule algebra, we have ${\cal A}\bowtie 
(\mf {A}\ot \mf {B})\simeq \mf{A} \gsl {\cal A}\trl \mf{B}$, hence, by 
Proposition \ref{crossed}, we obtain ${\cal A}\bowtie (\mf {A}\ot \mf {B})
\simeq {\cal A}\nat (\mf {A}\ot \mf {B})$. 
This raises the 
natural question whether we actually have ${\cal A}\bowtie {\mb A}\simeq 
{\cal A}\nat {\mb A}$ for any $H$-bimodule algebra ${\cal A}$ and any 
$H$-bicomodule algebra ${\mb A}$.   
We will see 
that this is the case, and for the proof we need first to recall  
some formulae from \cite{bpvo2}. Namely, if we denote by $\tilde{Q}^1_{\rho } 
\ot \tilde{Q}^2_{\rho }$ another copy of $\tilde{q}_{\rho }$, then we have: 
\begin{eqnarray}
&&
\Theta ^1_1\O ^1\ot \Theta ^1_2\O ^2\ot \tqra (\Theta ^2\O ^3)_{<0>}
\ot \O ^5\smi (\Theta ^3)_1(\tqrb )_1(\Theta ^2\O ^3)_{<1>_1}\nonumber \\
&&
\ot \;\O ^4\smi (\Theta ^3)_2(\tqrb )_2(\Theta ^2\O ^3)_{<1>_2}
=\tx ^1_{\l }\Theta ^1\ot \tx ^2_{\l }{\bf \Theta }^1  
\Theta ^2_{[-1]}\ov {\Theta }^1\ot 
\tx ^3_{\l }\tqra ({\bf \Theta }^2\Theta ^2_{[0]}\tQra 
\ov {\Theta }^2_{<0>})_{<0>}
\tx ^1_{\r }\nonumber \\
&&
\ot \;\smi ({\bf \Theta }^3\Theta ^3)
\tqrb ({\bf \Theta }^2\Theta ^2_{[0]}\tQra 
\ov {\Theta }^2_{<0>})_{<1>}\tx ^2_{\r }
\ot \smi (\ov {\Theta }^3)\tQrb \ov {\Theta }^2_{<1>}\tx ^3_{\r }\label{of2},
\end{eqnarray}  
\begin{eqnarray}
&&
\ov {\Theta }^1u_{<0>_{[-1]}}\ot (\tQra \ov {\Theta }^2_{<0>})_{<0>}
\tx ^1_{\r }u_{<0>_{[0]_{<0>}}}u^{'}_{<0>}\ot 
(\tQra \ov {\Theta }^2_{<0>})_{<1>}
\tx ^2_{\r }u_{<0>_{[0]_{<1>_1}}}u^{'}_{<1>_1}\nonumber\\
&&
\ot \;\smi (\ov {\Theta }^3u_{<1>})\tQrb \ov {\Theta }^2_{<1>}
\tx ^3_{\r }u_{<0>_{[0]_{<1>_2}}}u^{'}_{<1>_2}=
u_{[-1]}\ov {\Theta }^1\ot 
(u_{[0]}\tQra )_{<0>}(\ov {\Theta }^2u^{'})_{<0, 0>}\tx ^1_{\r }\nonumber\\
&&
\ot \;(u_{[0]}\tQra )_{<1>}(\ov {\Theta }^2u^{'})_{<0, 1>}\tx ^2_{\r }\ot   
\smi (\ov {\Theta }^3)\tQrb (\ov {\Theta }^2u^{'})_{<1>}\tx ^3_{\r }, 
\label{of3}
\end{eqnarray}
\begin{eqnarray}\label{of4} 
{\bf \Theta }^1\ot \tqra {\bf \Theta }^2_{<0>}\ot 
\smi ({\bf \Theta }^3)\tqrb {\bf \Theta }^2_{<1>}
=(\tqra )_{[-1]}
\theta ^1\ot (\tqra )_{[0]}\theta ^2\ot \tqrb \theta ^3.
\end{eqnarray} 
\begin{theorem} \label{isolrgd}
Let $H$ be a quasi-Hopf algebra, ${\cal A}$ an  
$H$-bimodule algebra and ${\mb A}$ an $H$-bicomodule algebra. Then the map 
\begin{eqnarray}
&&\nu :{\cal A}\bowtie {\mb A}\rightarrow {\cal A}\nat {\mb A},\\
&&\nu (\varphi \bowtie u)=\Theta ^1\cdot \varphi \cdot 
\smi (\Theta ^3)\tqrb \Theta ^2_{<1>}u_{<1>}\nat 
\tqra \Theta ^2_{<0>}u_{<0>}, \label{nu} 
\end{eqnarray}
for all $\varphi \in {\cal A}$ and $u\in {\mb A}$, is an algebra isomorphism, 
with inverse  
\begin{eqnarray}
&&\nu ^{-1}:{\cal A}\nat {\mb A}\rightarrow {\cal A}\bowtie {\mb A},\\
&&\nu ^{-1}(\varphi \nat u)=\theta ^1\cdot \varphi \cdot 
S^{-1}(\theta ^3u_{<1>}\tilde{p}^2_{\rho })\bowtie \theta ^2u_{<0>}
\tilde{p}^1_{\rho }, \label{inu}
\end{eqnarray}
and this isomorphism is compatible with the isomorphisms  
$(A\ot B)\nat {\mb A}\simeq A\gsm {\mb A}\gtl B$,  
$(A\ot B)\bowtie {\mb A}\simeq A\gsm {\mb A}\gtl B$ and 
${\cal A}\nat (\mf {A}\ot \mf {B})\simeq  
\mf {A}\gsl {\cal A}\trl \mf {B}$, ${\cal A}\bowtie (\mf {A}\ot \mf {B})
\simeq \mf {A}\gsl {\cal A}\trl \mf {B}$.
\end{theorem} 
\begin{proof}
First we establish that $\nu $ is an algebra map. We compute:\\[2mm] 
${\;\;\;}$
$\nu ((\varphi \bowtie u)(\psi \bowtie u^{'}))$
\begin{eqnarray*}
{\rm (\ref{gdp}, \ref{nu})}&=&(\Theta ^1_1\O ^1\cd \varphi \cdot 
\O ^5\smi (\Theta ^3)_1(\tqrb )_1
(\Theta ^2\O ^3)_{<1>_1}u_{<0>_{[0]_{<1>_1}}}u^{'}_{<1>_1})\\
&&(\Theta ^1_2\O ^2u_{<0>_{[-1]}}\cd \psi \cdot  
\smi (u_{<1>})\O ^4\smi (\Theta ^3)_2(\tqrb )_2
(\Theta ^2\O ^3)_{<1>_2}u_{<0>_{[0]_{<1>_2}}}u^{'}_{<1>_2})\\
&&\nat \tqra 
(\Theta ^2\O ^3)_{<0>}u_{<0>_{[0]_{<0>}}}u^{'}_{<0>}\\  
{\rm (\ref{of2})}&=&(\tx ^1_{\l }\Theta ^1\cd \varphi \cdot 
\smi ({\bf \Theta }^3\Theta ^3)\tqrb {\bf \Theta }^2_{<1>}
\Theta ^2_{[0]_{<1>}}(\tQra \ov {\Theta }^2_{<0>})_{<1>}\tx ^2_{\r }
u_{<0>_{[0]_{<1>_1}}}u^{'}_{<1>_1})\\
&&(\tx ^2_{\l }{\bf \Theta }^1\Theta ^2_{[-1]}\ov {\Theta }^1
u_{<0>_{[-1]}}\cd \psi \cdot 
\smi (\ov {\Theta }^3u_{<1>})\tQrb \ov {\Theta }^2_{<1>}\tx ^3_{\r }
u_{<0>_{[0]_{<1>_2}}}u^{'}_{<1>_2})\\
&&\nat \tx ^3_{\l }\tqra {\bf \Theta }^2_{<0>}
\Theta ^2_{[0]_{<0>}}(\tQra \ov {\Theta }^2_{<0>})_{<0>}
\tx ^1_{\r }u_{<0>_{[0]_{<0>}}}u^{'}_{<0>}\\
{\rm (\ref{of3})}&=&(\tx ^1_{\l }\Theta ^1\cd \varphi \cdot 
\smi ({\bf \Theta }^3\Theta ^3)\tqrb {\bf \Theta }^2_{<1>}
\Theta ^2_{[0]_{<1>}}
(u_{[0]}\tQra )_{<1>}(\ov {\Theta }^2u^{'})_{<0, 1>}\tx ^2_{\r })\\
&&(\tx ^2_{\l }
{\bf \Theta }^1\Theta ^2_{[-1]}u_{[-1]}\ov {\Theta }^1\cd \psi \cdot 
\smi (\ov {\Theta }^3)\tQrb (\ov {\Theta }^2u^{'})_{<1>}
\tx ^3_{\r })\\
&&\nat \tx ^3_{\l }\tqra {\bf \Theta }^2_{<0>}\Theta ^2_{[0]_{<0>}}
(u_{[0]}\tQra )_{<0>}(\ov {\Theta }^2u^{'})_{<0, 0>}\tx ^1_{\r }\\
{\rm (\ref{of4})}&=&(\tx ^1_{\l }\Theta ^1\cd \varphi \cdot 
\smi (\Theta ^3)\tqrb \theta ^3
\Theta ^2_{[0]_{<1>}}
(u_{[0]}\tQra )_{<1>}(\ov {\Theta }^2u^{'})_{<0, 1>}\tx ^2_{\r })\\
&&(\tx ^2_{\l }(\tilde{q}^1_{\rho })_{[-1]}\theta ^1
\Theta ^2_{[-1]}u_{[-1]}\ov {\Theta }^1\cd \psi \cdot 
\smi (\ov {\Theta }^3)\tQrb (\ov {\Theta }^2u^{'})_{<1>}
\tx ^3_{\r })\\
&&\nat \tx ^3_{\l }(\tqra )_{[0]}\theta ^2  
\Theta ^2_{[0]_{<0>}}
(u_{[0]}\tQra )_{<0>}(\ov {\Theta }^2u^{'})_{<0, 0>}\tx ^1_{\r }\\
{\rm (\ref{bca1})}&=&(\tx ^1_{\l }\Theta ^1\cd \varphi \cdot 
\smi (\Theta ^3)\tqrb \Theta ^2_{<1>}u_{<1>}\theta ^3 
(\tQra )_{<1>}(\ov {\Theta }^2u^{'})_{<0, 1>}\tx ^2_{\r })\\
&&(\tx ^2_{\l }(\tilde{q}^1_{\rho })_{[-1]}
\Theta ^2_{<0>_{[-1]}}u_{<0>_{[-1]}}
\theta ^1\ov {\Theta }^1\cd \psi \cdot 
\smi (\ov {\Theta }^3)\tQrb (\ov {\Theta }^2u^{'})_{<1>}
\tx ^3_{\r })\\
&&\nat \tx ^3_{\l }(\tqra )_{[0]}\Theta ^2_{<0>_{[0]}}u_{<0>_{[0]}}
\theta ^2  
(\tQra )_{<0>}(\ov {\Theta }^2u^{'})_{<0, 0>}\tx ^1_{\r }\\
{\rm (\ref{nat})}&=&(\Theta ^1\cdot \varphi \cdot  
S^{-1}(\Theta ^3)\tilde{q}^2_{\rho }\Theta ^2_{<1>}u_{<1>}\nat 
\tilde{q}^1_{\rho }\Theta ^2_{<0>}u_{<0>})\\
&&(\overline{\Theta }^1\cdot \psi \cdot S^{-1}(\overline{\Theta }^3)
\tilde{Q}^2_{\rho }\overline{\Theta }^2_{<1>}u'_{<1>}\nat 
\tilde{Q}^1_{\rho }\overline{\Theta }^2_{<0>}u'_{<0>})\\
{\rm (\ref{nu})}&=&\nu (\varphi \bowtie u)
\nu (\psi \bowtie u^{'}),     
\end{eqnarray*}
as needed. The fact that  
$\nu (1\bowtie 1)=1\nat 1$ is trivial.\\
We prove now that $\nu $ and $\nu ^{-1}$ are inverses. Indeed, we have:
\begin{eqnarray*}
\nu \nu ^{-1}(\varphi \nat u)&=&\Theta ^1\theta ^1\cdot \varphi \cdot 
S^{-1}(\theta ^3u_{<1>}\tilde{p}^2_{\rho })S^{-1}(\Theta ^3)
\tilde{q}^2_{\rho }\Theta ^2_{<1>}\theta ^2_{<1>}u_{<0>_{<1>}}
(\tilde{p}^1_{\rho })_{<1>}\\
&&\nat   
\tilde{q}^1_{\rho }\Theta ^2_{<0>}\theta ^2_{<0>}u_{<0>_{<0>}}
(\tilde{p}^1_{\rho })_{<0>}\\
{\rm (\ref{tpqr1a}, \ref{tpqr2a})}&=&\varphi \nat u, 
\end{eqnarray*}
\begin{eqnarray*}
\nu ^{-1}\nu \;(\varphi \bowtie u)&=&\theta ^1\Theta ^1\cdot \varphi \cdot   
S^{-1}(\Theta ^3)\tilde{q}^2_{\rho }\Theta ^2_{<1>}u_{<1>}
S^{-1}(\theta ^3(\tilde{q}^1_{\rho })_{<1>}\Theta ^2_{<0>_{<1>}}
u_{<0>_{<1>}}\tilde{p}^2_{\rho })\\
&&\bowtie \theta ^2(\tilde{q}^1_{\rho })_{<0>}\Theta ^2_{<0>_{<0>}}
u_{<0>_{<0>}}\tilde{p}^1_{\rho })\\
{\rm (\ref{tpqr1}, \ref{tpqr2})}&=&\theta ^1\Theta ^1\cdot \varphi 
\cdot S^{-1}(\theta ^3\Theta ^3)\bowtie \theta ^2\Theta ^2u\\
&=&\varphi \bowtie u,
\end{eqnarray*}
and we are done.
\end{proof}
\begin{examples}\rm 
${\;\;}$1) If $A$ is a left $H$-module algebra regarded as an    
$H$-bimodule algebra with trivial  
right $H$-action, then $A\bowtie {\mb A}$ and $A\nat {\mb A}$ both coincide 
with $A\gsm {\mb A}$, and the isomorphism $\nu $ is just the identity.\\ 
${\;\;}$2) If ${\mb A}=H$, the maps $\nu :{\cal A}\bowtie H\rightarrow 
{\cal A}\nat H$ and $\nu ^{-1}:{\cal A}\nat H\rightarrow  
{\cal A}\bowtie H$ are given by 
\begin{eqnarray*}
&&\nu (\varphi \bowtie h)=X^1\cdot \varphi \cdot  
\smi (X^3)q^2_RX^2_2h_2\nat  
q^1_RX^2_1h_1,  \\
&&\nu ^{-1}(\varphi \nat h)=x^1\cdot \varphi \cdot   
\smi (x^3h_2p^2_R)\bowtie x^2h_1p^1_R, 
\end{eqnarray*}
for all $\varphi \in {\cal A}$ and $h\in H$, where $q_R=q^1_R\ot q^2_R=
Y^1\ot S^{-1}(\alpha Y^3)Y^2$ and $p_R=p^1_R\ot p^2_R=y^1\ot y^2\beta 
S(y^3)$. \\  
${\;\;}$3) Let $H$ be a Hopf algebra with bijective antipode, ${\cal A}$ an 
$H$-bimodule algebra and ${\mb A}$ an $H$-bicomodule algebra in the usual 
(Hopf) sense. Then the maps 
$\nu :{\cal A}\bowtie {\mb A}\rightarrow  
{\cal A}\nat {\mb A}$ and $\nu ^{-1}:{\cal A}\nat {\mb A}\rightarrow   
{\cal A}\bowtie {\mb A}$ become:
\begin{eqnarray*}
&&\nu (\varphi \bowtie u)=\varphi \cdot u_{<1>}\nat u_{<0>},  \\
&&\nu ^{-1}(\varphi \nat u)=\varphi \cdot     
\smi (u_{<1>})\bowtie u_{<0>}, 
\end{eqnarray*}
for all $\varphi \in {\cal A}$ and $u\in {\mb A}$; if moreover ${\mb A}=H$,  
they become 
\begin{eqnarray*}
&&\nu (\varphi \bowtie h)=\varphi \cdot h_2\nat h_1,  \\
&&\nu ^{-1}(\varphi \nat h)=\varphi \cdot    
\smi (h_2)\bowtie h_1, 
\end{eqnarray*}
for all $\varphi \in {\cal A}$ and $h\in H$.
\end{examples}
Let now $H$ be a finite dimensional quasi-Hopf algebra. Recall that the  
quantum double $D(H)$ (generalizing the usual Drinfeld double of a 
Hopf algebra)  
was first introduced by Majid in \cite{m1} by an 
implicit Tannaka-Krein reconstruction procedure, and more explicit 
descriptions were obtained afterwards by Hausser and Nill in \cite{hn1},  
\cite{hn2}. According to one of these descriptions, the algebra 
structure of $D(H)$ is just the diagonal crossed product $H^*\bowtie H$.  
By transferring the whole structure of $D(H)$ via the map $\nu $, we can 
thus obtain a new realization of $D(H)$, having the L-R-smash product 
$H^*\nat H$ for the algebra structure. \\[2mm] 
We study the invariance under twisting of the L-R-smash product and  
first recall some facts from \cite{hn1}, \cite{bpvo2}. \\
Let $H$ be a quasi-bialgebra, ${\cal A}$  
an $H$-bimodule algebra and $F\in H\ot H$ a gauge transformation. 
If we introduce on ${\cal A}$  
another multiplication, by $\varphi \circ \varphi '=(G^1\cdot  
\varphi \cdot F^1)(G^2\cdot \varphi '\cdot F^2)$ for all $\varphi , 
\varphi '\in {\cal A}$, where $F^{-1}=G^1\ot G^2$,  
and denote this structure by $_{F}{\cal A}_{F^{-1}}$,  
then $_{F}{\cal A}_{F^{-1}}$ is an $H_F$-bimodule  
algebra, with the same unit and $H$-actions as ${\cal A}$.\\ 
Suppose that we have a left $H$-comodule algebra ${\mf B}$; then 
on the algebra structure of ${\mf B}$ one can introduce  
a left $H_F$-comodule algebra structure (denoted by  
${\mf B}^{F^{-1}}$ in what follows) putting  
$\lambda ^{F^{-1}}=\lambda $ and  
$\Phi _{\lambda }^{F^{-1}}=\Phi _{\lambda }(F^{-1}\ot 1_{\mf B})$. 
Similarly, if ${\mf A}$ is a right $H$-comodule algebra,  
one can introduce on the algebra structure of ${\mf A}$ a right 
$H_F$-comodule algebra structure (denoted by $^{F}{\mf A}$ in what follows)  
putting $^{F}\rho =\rho $ and $^{F}\Phi _{\rho }=(1_{\mf A}\ot F) 
\Phi _{\rho }$. 
One may check that if ${\mb A}$ is an $H$-bicomodule algebra, the  
left and right $H_F$-comodule algebras ${\mb A}^{F^{-1}}$ 
respectively $^{F}{\mb A}$ actually define the structure of an  
$H_F$-bicomodule algebra on ${\mb A}$, denoted by $^{F}{\mb A}^{F^{-1}}$, 
which has the same $\Phi _{\lambda , \rho }$ as ${\mb A}$. \\   
\begin{proposition}
With notation as above, we have an algebra isomorphism 
\begin{eqnarray*}
&&{\cal A}\nat {\mb A}\equiv \; _F{\cal A}_{F^{-1}}\nat ^F{\mb A}^{F^{-1}},
\end{eqnarray*}
given by the trivial identification. 
\end{proposition}
\begin{proof}
Let ${\cal F}^1\ot {\cal F}^2$ and ${\cal G}^1\ot {\cal G}^2$ be two more 
copies of $F$ and $F^{-1}$ respectively. We compute the multiplication 
in $_F{\cal A}_{F^{-1}}\nat ^F{\mb A}^{F^{-1}}$:\\[2mm]
${\;\;\;\;\;\;}$
$(\varphi \nat u)(\psi \nat u')$
\begin{eqnarray*}
&=&(F^1\tilde{x}^1_{\lambda }\cd \varphi \cd  
\theta ^3u'_{<1>}\tilde{x}^2_{\rho }G^1)\circ 
(F^2\tilde{x}^2_{\lambda }u_{[-1]}\theta ^1\cd \psi \cd  
\tilde{x}^3_{\rho }G^2)\nat \tilde{x}^3_{\lambda }u_{[0]}\theta ^2u'_{<0>} 
\tilde{x}^1_{\rho }\\
&=&({\cal G}^1F^1\tilde{x}^1_{\lambda }\cd \varphi \cd   
\theta ^3u'_{<1>}\tilde{x}^2_{\rho }G^1{\cal F}^1) 
({\cal G}^2F^2\tilde{x}^2_{\lambda }u_{[-1]}\theta ^1\cd \psi \cd   
\tilde{x}^3_{\rho }G^2{\cal F}^2)
\nat \tilde{x}^3_{\lambda }u_{[0]}\theta ^2u'_{<0>}\tilde{x}^1_{\rho }\\
&=&(\tilde{x}^1_{\lambda }\cd \varphi \cd   
\theta ^3u'_{<1>}\tilde{x}^2_{\rho }) 
(\tilde{x}^2_{\lambda }u_{[-1]}\theta ^1\cd \psi \cd   
\tilde{x}^3_{\rho })\nat \tilde{x}^3_{\lambda }u_{[0]}\theta ^2u'_{<0>} 
\tilde{x}^1_{\rho },
\end{eqnarray*}
which is the multiplication of ${\cal A}\nat {\mb A}$.
\end{proof}
\section{More properties of the L-R-smash product over 
bialgebras and Hopf algebras}
\setcounter{equation}{0}
The first property we want to emphasize is a direct consequence of 
Theorem \ref{isolrgd}.
\begin{proposition}
Let $H$ be a cocommutative Hopf algebra and ${\cal A}$ an $H$-bimodule 
algebra.  
Then the L-R-smash product ${\cal A}\nat H$ is isomorphic to an 
ordinary smash product ${\cal A}\# H$, where the left $H$-action on 
${\cal A}$ is now given by $h\rightarrow \varphi =h_1\cdot \varphi \cdot 
S(h_2)$, for all $h\in H$ and $\varphi \in {\cal A}$.
\end{proposition}
\begin{proof}
Since $H$ is cocommutative (hence in particular $S^2=id$), the generalized 
diagonal crossed product ${\cal A}\bowtie H$ (which is isomorphic to 
${\cal A}\nat H$) has multiplication:
\begin{eqnarray*}
(\varphi \bowtie h)(\varphi '\bowtie h')&=&
\varphi (h_1\cdot \varphi '\cdot S(h_3))\bowtie h_2h'\\
&=&\varphi (h_1\cdot \varphi '\cdot S(h_2))\bowtie h_3h'\\
&=&\varphi (h_1\rightarrow \varphi ')\bowtie h_2h',
\end{eqnarray*}
which is just the multiplication of ${\cal A}\# H$.
\end{proof}
The following provides us with a Maschke-type theorem 
for L-R-smash products, the proof is  
similar to the one for crossed products, see \cite{mon}, Theorem 7.4.2.
\begin{theorem}\label{maschke}
Let $H$ be a finite dimensional Hopf algebra such that $H$ is semisimple 
and $H^*$ is unimodular, and let ${\cal A}$ be an $H$-bimodule algebra. 
Then:\\
(i) If $V\in \;_{{\cal A}\nat H}{\cal M}$ and $W\subseteq V$ is a   
submodule which has a complement in $_{\cal A}{\cal M}$, then $W$ has   
also a complement in $\;_{{\cal A}\nat H}{\cal M}$.\\
(ii) If ${\cal A}$ is semisimple Artinian, then so is ${\cal A}\nat H$.
\end{theorem}
\begin{proof}
Obviously (ii) follows from (i), so we prove (i) now. Let $t\in H$ be an 
integral with $\varepsilon (t)=1$. By \cite{radford}, formula (15), we 
have 
\begin{eqnarray*}
&&S^{-1}(t_3)g^{-1}t_1\ot t_2=1\ot t, 
\end{eqnarray*}
where $g$ is the distinguished group-like element of $H$; our 
hypothesis that $H^*$ is unimodular implies $g=1$, hence we also get     
\begin{eqnarray}
&&S^{-1}(t_4)t_1\ot t_2\ot t_3=1\ot t_1\ot t_2. \label{uni} 
\end{eqnarray}
Let now $\pi :V\rightarrow W$ be an ${\cal A}$-projection. We construct the 
averaging function as in \cite{mon}, namely:
\begin{eqnarray*}
&&\tilde{\pi }(v)=(1\nat S(t_1))\cdot \pi ((1\nat t_2)\cdot v), \;\;\;
\forall \;v\in V.
\end{eqnarray*}
Before proving that $\tilde{\pi }$ is ${\cal A}$-linear, note the 
following formulae:
\begin{eqnarray}
&&(1\nat h)(\varphi \nat 1)=(h_1\cdot \varphi \cdot S^{-1}(h_3)\nat 1)
(1\nat h_2), \label{int1} \\
&&\varphi \nat h=(\varphi \cdot S^{-1}(h_2)\nat 1)(1\nat h_1), \label{int2}
\end{eqnarray}
for all $\varphi \in {\cal A}$ and $h\in H$, which follow by direct 
computation using the formula for the multiplication in  
${\cal A}\nat H$. 
Let now $\varphi \in {\cal A}$ and $v\in V$; we compute:
\begin{eqnarray*}
\tilde{\pi }((\varphi \nat 1)\cdot v)&=&(1\nat S(t_1))\cdot 
\pi ((1\nat t_2)(\varphi \nat 1)\cdot v)\\
{\rm (\ref{int1})}&=&(1\nat S(t_1))\cdot \pi ((t_2\cdot \varphi \cdot 
S^{-1}(t_4)\nat 1)(1\nat t_3)\cdot v)\\
&=&(1\nat S(t_1))(t_2\cdot \varphi \cdot S^{-1}(t_4)\nat 1)\cdot 
\pi ((1\nat t_3)\cdot v)\;\;\;(since\;\pi \;is\;{\cal A}-linear)\\
&=&(S(t_1)_1t_2\cdot \varphi \cdot S^{-1}(t_4)\nat S(t_1)_2)\cdot 
\pi ((1\nat t_3)\cdot v)\\
&=&(\varphi \cdot S^{-1}(t_3)\nat S(t_1))\cdot \pi ((1\nat t_2)\cdot v)\\
{\rm (\ref{int2})}&=&(\varphi \cdot S^{-1}(t_3)S^{-1}(S(t_1)_2)\nat 1)
(1\nat S(t_1)_1)\cdot \pi ((1\nat t_2)\cdot v)\\
&=&(\varphi \cdot S^{-1}(t_4)t_1\nat 1)(1\nat S(t_2))\cdot \pi ((1\nat t_3)
\cdot v)\\
{\rm (\ref{uni})}&=&(\varphi \nat 1)(1\nat S(t_1))\cdot \pi ((1\nat t_2)
\cdot v)\\
&=&(\varphi \nat 1)\cdot \tilde{\pi }(v),
\end{eqnarray*}
hence $\tilde{\pi }$ is ${\cal A}$-linear. The rest of the proof is 
identical to the one in \cite{mon}, Theorem 7.4.2.
\end{proof}
As a consequence of Theorem \ref{isolrgd} and Theorem \ref{maschke}, 
we obtain the very well-known result (see \cite{radford}, Proposition 7):
\begin{corollary}
If $H$ is a semisimple cosemisimple Hopf algebra, then $D(H)$ is semisimple. 
\end{corollary}
Let now $H$ be a bialgebra, ${\cal A}$ an $H$-bimodule algebra and 
${\mb A}$ an $H$-bicomodule algebra in the usual (Hopf) sense. Let also 
$A$ be an algebra in the Yetter-Drinfeld category $_H^H{\cal YD}$, that is 
$A$ is a left $H$-module algebra, a left $H$-comodule algebra (with left 
$H$-comodule structure denoted by $a\mapsto a_{(-1)}\ot a_{(0)}\in H\ot A$) 
and the Yetter-Drinfeld compatibility condition holds: 
\begin{eqnarray}
&&h_1a_{(-1)}\ot h_2\cdot a_{(0)}=(h_1\cdot a)_{(-1)}h_2\ot 
(h_1\cdot a)_{(0)}, \;\;\;\forall \;h\in H, \;a\in A. \label{yd}
\end{eqnarray}
Consider first the generalized smash product ${\cal A}\gsm A$, an  
associative algebra. From the condition (\ref{yd}), it follows   
that ${\cal A}\gsm A$ becomes an $H$-bimodule algebra,  
with $H$-actions
\begin{eqnarray*}
&&h\cdot (\varphi \gsm a)=h_1\cdot \varphi \gsm h_2\cdot a, \\
&&(\varphi \gsm a)\cdot h=\varphi \cdot h\gsm a,
\end{eqnarray*}
for all $h\in H$, $\varphi \in {\cal A}$ and $a\in A$, hence we may consider  
the algebra $({\cal A}\gsm A)\nat {\mb A}$.\\
Then, consider the generalized smash product $A\gsm {\mb A}$, an  
associative algebra. Using the condition (\ref{yd}), one can see that  
$A\gsm {\mb A}$ becomes an $H$-bicomodule algebra, with $H$-coactions
\begin{eqnarray*}
&&\rho :A\gsm {\mb A}\rightarrow (A\gsm {\mb A})\ot H,\;\;\;\rho (a\gsm u)=
(a\gsm u_{<0>})\ot u_{<1>}, \\
&&\lambda :A\gsm {\mb A}\rightarrow H\ot (A\gsm {\mb A}), \;\;\;
\lambda (a\gsm u)=a_{(-1)}u_{[-1]}\ot (a_{(0)}\gsm u_{[0]}), 
\end{eqnarray*}
for all $a\in A$ and $u\in {\mb A}$, hence we may consider the algebra  
${\cal A}\nat (A\gsm {\mb A})$. 
\begin{proposition}
We have an algebra isomorphism $({\cal A}\gsm A)\nat {\mb A}\equiv 
{\cal A}\nat (A\gsm {\mb A})$, given by the trivial identification. 
\end{proposition}
\begin{proof}
The multiplication in $({\cal A}\gsm A)\nat {\mb A}$ is:\\[2mm]
${\;\;\;\;\;\;\;\;\;\;\;}$
$((\varphi \gsm a)\nat u)((\varphi '\gsm a')\nat u')$
\begin{eqnarray*}
&&=((\varphi \gsm a)\cdot u'_{<1>})(u_{[-1]}\cdot (\varphi '\gsm a'))\nat 
u_{[0]}u'_{<0>}\\
&&=(\varphi \cdot u'_{<1>}\gsm a)(u_{[-1]_1}\cdot \varphi '\gsm 
u_{[-1]_2}\cdot a')\nat u_{[0]}u'_{<0>}\\
&&=((\varphi \cdot u'_{<1>})(a_{(-1)}u_{[-1]_1}\cdot \varphi ')\gsm a_{(0)}
(u_{[-1]_2}\cdot a'))\nat u_{[0]}u'_{<0>}.
\end{eqnarray*}
The multiplication in ${\cal A}\nat (A\gsm {\mb A})$ is:\\[2mm]
${\;\;\;\;\;\;\;\;\;\;\;}$
$(\varphi \nat (a\gsm u))(\varphi '\nat (a'\gsm u'))$
\begin{eqnarray*}
&&=(\varphi \cdot (a'\gsm u')_{<1>})((a\gsm u)_{[-1]}\cdot \varphi ')\nat 
(a\gsm u)_{[0]}(a'\gsm u')_{<0>}\\
&&=(\varphi \cdot u'_{<1>})(a_{(-1)}u_{[-1]}\cdot \varphi ')\nat 
(a_{(0)}\gsm u_{[0]})(a'\gsm u'_{<0>})\\
&&=(\varphi \cdot u'_{<1>})(a_{(-1)}u_{[-1]}\cdot \varphi ')\nat 
(a_{(0)}(u_{[0]_{[-1]}}\cdot a')\gsm u_{[0]_{[0]}}u'_{<0>})\\
&&=(\varphi \cdot u'_{<1>})(a_{(-1)}u_{[-1]_1}\cdot \varphi ')\nat  
(a_{(0)}(u_{[-1]_2}\cdot a')\gsm u_{[0]}u'_{<0>}), 
\end{eqnarray*}
hence the two multiplications coincide.
\end{proof}
Since the L-R-smash product coincides with the generalized smash product 
if the right $H$-action is trivial, we also obtain:
\begin{corollary}
If $H$, $A$, ${\mb A}$ are as above and $A'$ is a left $H$-module algebra, 
then we have an algebra isomorphism $(A'\gsm A)\gsm {\mb A}\equiv 
A'\gsm (A\gsm {\mb A})$, given by the trivial identification. 
\end{corollary}
\section{L-R-twisting data}
\setcounter{equation}{0}
We start by recalling the set-up used in \cite{fst}, slightly  
modifying terminology.\\
Let $(A, \mu )$ be a (unital) associative algebra. Assume that there exists 
a bialgebra $H$ such that $A$ is a left $H$-module algebra (denote by 
$\pi :H\ot A\rightarrow A$, $\pi (h\ot a)=h\cdot a$ the action), a left 
$H$-comodule algebra (denote by $\psi :A\rightarrow H\ot A$, 
$\psi (a)=a_{(-1)}\ot a_{(0)}$ the coaction) and the following 
compatibility condition holds:
\begin{eqnarray}
&&(h\cdot a)_{(-1)}\ot (h\cdot a)_{(0)}=a_{(-1)}\ot h\cdot a_{(0)}, \;\;\;
\forall \;\;h\in H,\;a\in A. \label{long}
\end{eqnarray}
We call the triple $(H, \pi , \psi )$ a left twisting datum for $(A, \mu )$ 
(in \cite{fst} it is called a very strong left twisting datum). If we 
define a new multiplication on $A$, by 
\begin{eqnarray}
&&a\star b=a_{(0)}(a_{(-1)}\cdot b), \;\;\;\forall 
\;\;a, b\in A, \label{star}
\end{eqnarray}
then this multiplication defines a new algebra structure on $A$, 
with the same unit. The product $\star $ is called the left twisted 
product.  
\begin{example}{\em \label{ex}
Let $H$ be a Hopf algebra with bijective antipode $S$, ${\cal A}$ an 
$H$-bimodule algebra with actions $h\ot \varphi \mapsto h\cdot \varphi $ 
and $\varphi \ot h\mapsto \varphi \cdot h$ for all $h\in H$, 
$\varphi \in {\cal A}$, and ${\mb A}$ an $H$-bicomodule algebra with 
coactions $u\mapsto u_{[-1]}\ot u_{[0]}$, $u\mapsto u_{<0>}\ot 
u_{<1>}$ for all $u\in {\mb A}$, and denote also by $u_{\{-1\}}\ot 
u_{\{0\}}\ot u_{\{1\}}:=u_{<0>_{[-1]}}\ot u_{<0>_{[0]}}\ot u_{<1>}=
u_{[-1]}\ot u_{[0]_{<0>}}\ot u_{[0]_{<1>}}$. Then ${\cal A}$ becomes a 
left $H\ot H^{op}$-module algebra with action $\pi :H\ot H^{op}\ot 
{\cal A}\rightarrow {\cal A}$, $\pi (h\ot h'\ot \varphi )=h\cdot \varphi  
\cdot h'$, and ${\mb A}$ becomes a left $H\ot H^{op}$-comodule algebra 
(see \cite{stef}) with coaction $\psi :{\mb A}\rightarrow (H\ot H^{op})\ot 
{\mb A}$, $u\mapsto (u_{\{-1\}}\ot S^{-1}(u_{\{1\}}))\ot u_{\{0\}}$. 
It is easy to check that $(H\ot H^{op}, \pi \ot id, (\tau \ot id)(id \ot 
\psi ))$, where $\tau $ is the usual twist, is a left twisting datum 
for ${\cal A}\ot {\mb A}$, and the corresponding twisted product is 
\begin{eqnarray*}
&&(\varphi \ot u)\star (\varphi '\ot u')=\varphi (u_{\{-1\}}\cdot 
\varphi '\cdot S^{-1}(u_{\{1\}}))\ot u_{\{0\}}u', 
\end{eqnarray*}
and this is exactly the multiplication of the generalized diagonal 
crossed product ${\cal A}\bowtie {\mb A}$.}
\end{example}  
Now we introduce the L-R-version of the above construction. Let $(A, \mu )$  
be an algebra and assume there exists a bialgebra $H$ such that:\\
(i) $A$ is an $H$-bimodule algebra with actions denoted by 
$\pi _l:H\ot A\rightarrow A$, $\pi _l(h\ot a)=h\cdot a$ and 
$\pi _r:A\ot H\rightarrow A$, $\pi _r(a\ot h)=a\cdot h$; \\
(ii) $A$ is an $H$-bicomodule algebra, with coactions denoted by 
$\psi _l:A\rightarrow H\ot A$, $a\mapsto a_{[-1]}\ot a_{[0]}$ and 
$\psi _r:A\rightarrow A\ot H$, $a\mapsto a_{<0>}\ot a_{<1>}$; \\
(iii) The following compatibility conditions hold, for all $h\in H$ and 
$a\in A$:
\begin{eqnarray}
&&(h\cdot a)_{[-1]}\ot (h\cdot a)_{[0]}=a_{[-1]}\ot h\cdot a_{[0]}, 
\label{a1}\\
&&(h\cdot a)_{<0>}\ot (h\cdot a)_{<1>}=h\cdot a_{<0>}\ot a_{<1>},  
\label{a2}\\
&&(a\cdot h)_{[-1]}\ot (a\cdot h)_{[0]}=a_{[-1]}\ot a_{[0]}\cdot h, 
\label{a3}\\
&&(a\cdot h)_{<0>}\ot (a\cdot h)_{<1>}=a_{<0>}\cdot h\ot a_{<1>}. 
\label{a4}
\end{eqnarray}
We call $(H, \pi _l, \pi _r, \psi _l, \psi _r)$ an L-R-twisting datum for 
$A$. Given such a datum, we define a new multiplication on $A$, by  
\begin{eqnarray}
&&a\bullet b=(a_{[0]}\cdot b_{<1>})(a_{[-1]}\cdot b_{<0>}), \;\;\;
\forall \;\;a, b\in A, \label{bullet}
\end{eqnarray}
and call it the L-R-twisted product. Obviously it has the same unit $1$ 
as $A$. Note the following  
easy consequences of the axioms:
\begin{eqnarray}
&&h\cdot (a\bullet b)=(h_1\cdot a_{[0]}\cdot b_{<1>})(h_2a_{[-1]}\cdot 
b_{<0>}), \label{c1} \\
&&(a\bullet b)\cdot h=(a_{[0]}\cdot b_{<1>}h_1)(a_{[-1]}\cdot b_{<0>} 
\cdot h_2), \label{c2} \\
&&(a\bullet b)_{[-1]}\ot (a\bullet b)_{[0]}=a_{[0]_{[-1]}}
b_{<0>_{[-1]}}\ot (a_{[0]_{[0]}}\cdot b_{<1>})(a_{[-1]}\cdot 
b_{<0>_{[0]}}), \label{c3} \\
&&(a\bullet b)_{<0>}\ot (a\bullet b)_{<1>}=(a_{[0]_{<0>}}\cdot b_{<1>})
(a_{[-1]}\cdot b_{<0>_{<0>}})\ot a_{[0]_{<1>}}b_{<0>_{<1>}}. 
\label{c4}
\end{eqnarray}  
\begin{proposition}
$(A, \bullet , 1)$ is an associative unital algebra.
\end{proposition}
\begin{proof}
We compute:
\begin{eqnarray*}
(a\bullet b)\bullet c&=&((a\bullet b)_{[0]}\cdot c_{<1>})
((a\bullet b)_{[-1]}\cdot c_{<0>})\\
{\rm (\ref{c3})}&=&(a_{[0]_{[0]}}\cdot b_{<1>}c_{<1>_1})
(a_{[-1]}\cdot b_{<0>_{[0]}}\cdot c_{<1>_2})(a_{[0]_{[-1]}}b_{<0>_{[-1]}}
\cdot c_{<0>}), 
\end{eqnarray*}
\begin{eqnarray*}
a\bullet (b\bullet c)&=&(a_{[0]}\cdot (b\bullet c)_{<1>})
(a_{[-1]}\cdot (b\bullet c)_{<0>})\\
{\rm (\ref{c4})}&=&(a_{[0]}\cdot b_{[0]_{<1>}}c_{<0>_{<1>}})
(a_{[-1]_1}\cdot b_{[0]_{<0>}}\cdot c_{<1>})
(a_{[-1]_2}b_{[-1]}\cdot c_{<0>_{<0>}}), 
\end{eqnarray*}
and the two terms are equal because $A$ is an $H$-bicomodule.
\end{proof}
\begin{remark}{\em 
An L-R-twisting datum is in particular a left twisting datum, but in general 
the corresponding twisted products $\bullet $ and respectively $\star $  
are different. On the other hand, any left twisting datum can be regarded 
as an L-R-twisting datum with trivial right action and coaction, 
and in this case the corresponding twisted products coincide.} 
\end{remark}
As a particular case of an L-R-twisting datum, obtained if the left 
action and coaction are trivial, we obtain the following concept. 
\begin{definition}
Let $(A, \mu )$ be a (unital) associative algebra. Assume that there exists 
a bialgebra $H$ such that $A$ is a right $H$-module algebra (denote by  
$\pi _r:A\ot H\rightarrow A$, $\pi _r(a\ot h)=a\cdot h$ the action), a right  
$H$-comodule algebra (denote by $\psi _r:A\rightarrow A\ot H$,  
$\psi _r(a)=a_{<0>}\ot a_{<1>}$ the coaction) and the following  
compatibility condition holds, for all $a\in A$ and $h\in H$:
\begin{eqnarray*}
&&(a\cdot h)_{<0>}\ot (a\cdot h)_{<1>}=a_{<0>}\cdot h\ot a_{<1>}. 
\end{eqnarray*}
We call the triple $(H, \pi _r, \psi _r)$ a right twisting datum 
for $(A, \mu )$. If we 
define a new multiplication on $A$, by 
\begin{eqnarray}
&&a\diamond b=(a\cdot b_{<1>})b_{<0>}, \;\;\;\forall 
\;\;a, b\in A, \label{diamond}
\end{eqnarray}
then this multiplication defines a new algebra structure on $A$, 
with the same unit. The product $\diamond $ is called the right twisted  
product.  
\end{definition}
\begin{proposition}\label{iterate}
Let $A$ be an algebra and $(H, \pi _l, \pi _r, \psi _l, \psi _r)$ an 
L-R-twisting datum for $A$, with notation as before. Then the L-R-twisted 
product $\bullet $ can be obtained as a left twisting followed 
by a right twisting and also viceversa. 
\end{proposition}
\begin{proof}
First consider the left twisted product algebra $(A, \star )$; it is easy 
to see that $(H, \pi _r, \psi _r)$ is a right twisting datum for 
$(A, \star )$, and the corresponding right twisted product becomes:
\begin{eqnarray*}
a\diamond b&=&(a\cdot b_{<1>})\star b_{<0>}\\
&=&(a\cdot b_{<1>})_{[0]}((a\cdot b_{<1>})_{[-1]}\cdot b_{<0>})\\
{\rm (\ref{a3})}&=&(a_{[0]}\cdot b_{<1>})(a_{[-1]}\cdot b_{<0>})\\
&=&a\bullet b.
\end{eqnarray*}
Similarly, one can start with the right twisted product algebra 
$(A, \diamond )$, for which $(H, \pi _l, \psi _l)$ is a left twisting  
datum,  and the corresponding left twisted product 
coincides with the L-R-twisted product.
\end{proof}
\begin{example}{\em 
Let $H$ be a bialgebra, 
${\cal A}$ an 
$H$-bimodule algebra with actions $h\ot \varphi \mapsto h\cdot \varphi $ 
and $\varphi \ot h\mapsto \varphi \cdot h$ for all $h\in H$, 
$\varphi \in {\cal A}$, and ${\mb A}$ an $H$-bicomodule algebra with 
coactions $u\mapsto u_{[-1]}\ot u_{[0]}$, $u\mapsto u_{<0>}\ot 
u_{<1>}$ for all $u\in {\mb A}$. Take the algebra $A={\cal A}\ot 
{\mb A}$, which becomes an $H$-bimodule algebra with actions 
$h\cdot (\varphi \ot u)=h\cdot \varphi \ot u$ and 
$(\varphi \ot u)\cdot h=\varphi \cdot h\ot u$, for all $h\in H$, 
$\varphi \in {\cal A}$, $u\in {\mb A}$, and an $H$-bicomodule algebra, 
with coactions $\varphi \ot u\mapsto u_{[-1]}\ot (\varphi \ot u_{[0]})$, 
$\varphi \ot u\mapsto (\varphi \ot u_{<0>})\ot u_{<1>}$. Moreover, one  
checks that the conditions (\ref{a1})-(\ref{a4}) are satisfied,  
hence we have an L-R-twisting datum for ${\cal A}\ot {\mb A}$. The 
corresponding L-R-twisted product is: 
\begin{eqnarray*}
(\varphi \ot u)\bullet (\varphi '\ot u')&=&
((\varphi \ot u)_{[0]}\cdot (\varphi '\ot u')_{<1>})
((\varphi \ot u)_{[-1]}\cdot (\varphi '\ot u')_{<0>})\\
&=&((\varphi \ot u_{[0]})\cdot u'_{<1>})(u_{[-1]}\cdot (\varphi '
\ot u'_{<0>}))\\
&=&(\varphi \cdot u'_{<1>})(u_{[-1]}\cdot \varphi ')\ot u_{[0]}u'_{<0>}, 
\end{eqnarray*}
and this is exactly the multiplication of the L-R-smash product 
${\cal A}\nat {\mb A}$.}
\end{example}  
If $H$ is a Hopf algebra with bijective antipode $S$, ${\cal A}$ is an 
$H$-bimodule algebra and ${\mb A}$ an $H$-bicomodule algebra, 
we have proved in Section \ref{sectiune} that ${\cal A}\nat {\mb A}\simeq  
{\cal A}\bowtie {\mb A}$ as algebras. We derive now two    
interpretations of this result at the level of twisting data and  
twisted products. 
\begin{theorem}\label{iso1}
With notation as above, let $A$ be an algebra and 
$(H, \pi _l, \pi _r, \psi _l, \psi _r)$ an  
L-R-twisting datum for $A$, $H$ being a    
Hopf algebra with bijective antipode $S$. If we denote by 
$\pi $ respectively  
$\psi $ the left $H\ot H^{op}$-module algebra respectively comodule  
algebra structures on $A$ defined as in Example \ref{ex}, then 
$(H\ot H^{op}, \pi , \psi )$ is a left twisting datum for $A$. 
Moreover, the corresponding twisted algebras $(A, \bullet )$ and 
$(A, \star )$ are isomorphic, and the isomorphism is defined by: 
\begin{eqnarray*}
&&\lambda :(A, \bullet )\rightarrow (A, \star ), \;\;
\lambda (a)=a_{<0>}\cdot S^{-1}(a_{<1>}), \;\;\;\forall \;\;a\in A, \\
&&\lambda ^{-1}:(A, \star )\rightarrow (A, \bullet ), \;\;
\lambda ^{-1}(a)=a_{<0>}\cdot a_{<1>}, \;\;\;\forall \;\;a\in A.
\end{eqnarray*} 
In particular, for $A={\cal A}\ot {\mb A}$, we obtain 
${\cal A}\nat {\mb A}\simeq {\cal A}\bowtie {\mb A}$. 
\end{theorem}  
\begin{proof}
To prove that $(H\ot H^{op}, \pi , \psi )$ is a left twisting datum for 
$A$, one has to check \ref{long}, and this follows using 
(\ref{a1})-(\ref{a4}). We only prove that $\lambda $ is an algebra 
isomorphism. The fact that $\lambda \lambda ^{-1}=\lambda ^{-1}\lambda =id$ 
follows easily using (\ref{a4}); obviously $\lambda (1)=1$, hence we only 
have to check that $\lambda $ is multiplicative:
\begin{eqnarray*}
\lambda (a\bullet b)&=&(a\bullet b)_{<0>}\cdot S^{-1}((a\bullet b)_{<1>})\\
{\rm (\ref{c4})}&=&[(a_{[0]_{<0>}}\cdot b_{<1>})(a_{[-1]}\cdot 
b_{<0>_{<0>}})]\cdot S^{-1}(a_{[0]_{<1>}}b_{<0>_{<1>}})\\
&=&[a_{[0]_{<0>}}\cdot b_{<1>}S^{-1}(b_{<0>_{<1>_2}})S^{-1}
(a_{[0]_{<1>_2}})]\\
&&[a_{[-1]}\cdot b_{<0>_{<0>}}\cdot S^{-1}(b_{<0>_{<1>_1}})
S^{-1}(a_{[0]_{<1>_1}})]\\
&=&[a_{[0]_{<0>}}\cdot b_{<1>_3}S^{-1}(b_{<1>_2})S^{-1}
(a_{[0]_{<1>_2}})]\\
&&[a_{[-1]}\cdot b_{<0>}\cdot S^{-1}(b_{<1>_1})
S^{-1}(a_{[0]_{<1>_1}})]\\
&=&[a_{[0]_{<0>}}\cdot S^{-1}(a_{[0]_{<1>_2}})]
[a_{[-1]}\cdot b_{<0>}\cdot S^{-1}(b_{<1>})S^{-1}(a_{[0]_{<1>_1}})],
\end{eqnarray*}
\begin{eqnarray*}
\lambda (a)\star \lambda (b)&=&
[a_{<0>}\cdot S^{-1}(a_{<1>})]\star [b_{<0>}\cdot S^{-1}(b_{<1>})]\\
&=&(a_{<0>}\cdot S^{-1}(a_{<1>}))_{<0>_{[0]}}
[(a_{<0>}\cdot S^{-1}(a_{<1>}))_{<0>_{[-1]}}\cdot b_{<0>}\cdot \\
&&\cdot S^{-1}(b_{<1>})S^{-1}((a_{<0>}\cdot S^{-1}(a_{<1>}))_{<1>})]\\
{\rm (\ref{a4})}&=&(a_{<0>_{<0>}}\cdot S^{-1}(a_{<1>}))_{[0]}
[(a_{<0>_{<0>}}\cdot S^{-1}(a_{<1>}))_{[-1]}\cdot b_{<0>}\cdot \\
&&\cdot S^{-1}(b_{<1>})S^{-1}(a_{<0>_{<1>}})]\\
{\rm (\ref{a3})}&=&[a_{<0>_{<0>_{[0]}}}\cdot S^{-1}(a_{<1>})]
[a_{<0>_{<0>_{[-1]}}}\cdot b_{<0>}\cdot S^{-1}(b_{<1>})
S^{-1}(a_{<0>_{<1>}})]\\
&=&[a_{<0>_{[0]}}\cdot S^{-1}(a_{<1>_2})][a_{<0>_{[-1]}}\cdot b_{<0>}\cdot  
S^{-1}(b_{<1>})S^{-1}(a_{<1>_1})]\\
&=&[a_{[0]_{<0>}}\cdot S^{-1}(a_{[0]_{<1>_2}})]
[a_{[-1]}\cdot b_{<0>}\cdot S^{-1}(b_{<1>})S^{-1}(a_{[0]_{<1>_1}})],
\end{eqnarray*}
and the proof is finished.
\end{proof}
The second interpretation is inspired by the following two facts. 
First, by Proposition \ref{iterate}, one can see that ${\cal A}\nat {\mb A}$ 
may be written as a right twisting of the generalized smash product 
${\cal A}\gsm {\mb A}$. Second, we have the observation in 
\cite{fs} that the Drinfeld double can be obtained as a two step 
twisting procedure, where at the first step the smash product is obtained; 
this can be extended to a generalized diagonal crossed product 
${\cal A}\bowtie {\mb A}$, thus obtaining ${\cal A}\bowtie {\mb A}$ 
as a left twisting of ${\cal A}\gsm {\mb A}$. Hence, our general result 
looks as follows:
\begin{theorem}
Let $A$ be an algebra and $(H, \pi _r, \psi _r)$ a right twisting datum 
for $A$, with notation as above, where $H$ is a Hopf algebra with 
bijective antipode. Define $\pi _l:H^{op}\ot A\rightarrow A$, 
$\pi _l(h\ot a)=a\cdot h$ and $\psi _l:A\rightarrow H^{op}\ot A$, 
$a\mapsto S^{-1}(a_{<1>})\ot a_{<0>}$. Then $(H^{op}, \pi _l, \psi _l)$ 
is a left twisting datum for $A$, and the corresponding twisted algebras  
$(A, \diamond )$ and $(A, \star )$ are isomorphic, via the maps: 
\begin{eqnarray*}
&&\lambda :(A, \diamond )\rightarrow (A, \star ), \;\;
\lambda (a)=a_{<0>}\cdot S^{-1}(a_{<1>}), \;\;\;\forall \;\;a\in A, \\
&&\lambda ^{-1}:(A, \star )\rightarrow (A, \diamond ), \;\;
\lambda ^{-1}(a)=a_{<0>}\cdot a_{<1>}, \;\;\;\forall \;\;a\in A.
\end{eqnarray*} 
In particular, for $A={\cal A}\gsm {\mb A}$, we obtain  
${\cal A}\nat {\mb A}\simeq {\cal A}\bowtie {\mb A}$. 
\end{theorem}  
\begin{proof}
Similar to the proof of Theorem \ref{iso1}.
\end{proof}
We end this section with a partial answer to the following natural question. 
Suppose that $A$ is an algebra and $(H, \pi , \psi )$ is a left twisting 
datum for $A$; then how far is the twisted algebra $(A, \star )$ from being 
isomorphic to a generalized smash product?
\begin{proposition}
Let $A$ and $(H, \pi , \psi )$ be as above. Define the algebras 
$A^H=\{a\in A/h\cdot a=\varepsilon (h)a,\;\forall \;h\in H\}$ and  
$A^{co (H)}=\{a\in A/a_{(-1)}\ot a_{(0)}=1\ot a\}$. Then $A^H$ is a left 
$H$-comodule algebra and $A^{co (H)}$ is a left $H$-module algebra. 
If moreover we have that $ab=ba$ for all $a\in A^{co (H)}$ and 
$b\in A^H$, then the map $\lambda :A^{co (H)}\gsm A^H\rightarrow 
(A, \star )$, $\lambda (a\ot b)=ab$, is an algebra map. 
\end{proposition} 
\begin{proof}
The assertions concerning $A^H$ and $A^{co (H)}$ follow easily from 
(\ref{long}). Assume now that 
\begin{eqnarray}
&&ab=ba, \;\;\forall \;\;a\in A^{co (H)}, \;b\in A^H. \label{comut}
\end{eqnarray} 
Then we compute, for all $a, a'\in A^{co (H)}$ and $b, b'\in A^H$:
\begin{eqnarray*}
\lambda (a\gsm b)\star \lambda (a'\gsm b')&=&(ab)\star (a'b')\\ 
&=&a_{(0)}b_{(0)}(a_{(-1)}b_{(-1)}\cdot (a'b'))\\
&=&ab_{(0)}(b_{(-1)}\cdot (a'b'))\;\;\;\;(since \;a\in A^{co (H)})\\
&=&ab_{(0)}(b_{(-1)_1}\cdot a')(b_{(-1)_2}\cdot b')\\
&=&ab_{(0)}(b_{(-1)}\cdot a')b'\;\;\;\;(since \;b'\in A^H)\\
&=&a(b_{(-1)}\cdot a')b_{(0)}b'\;\;\;\;(by \;(\ref{comut}))\\
&=&\lambda (a(b_{(-1)}\cdot a')\gsm b_{(0)}b')\\
&=&\lambda ((a\gsm b)(a'\gsm b')), 
\end{eqnarray*}
hence $\lambda $ is an algebra map.
\end{proof}
\section{L-R-smash coproduct over bialgebras}
\setcounter{equation}{0}
Throughout this section, $H$ will be a given bialgebra. We introduce  
the L-R-smash coproduct $C\nat H$, dualizing the L-R-smash product 
${\cal A}\nat H$ and generalizing Molnar's smash coproduct. \\[2mm]
Let $C$ be an $H$-bicomodule coalgebra, that is: \\
(i) $C$ is an $H$-bicomodule, with structures
\begin{eqnarray*}
&&\rho :C\rightarrow C\ot H,\;\;\rho (c)=c^{<0>}\ot c^{<1>}, \\
&&\lambda :C\rightarrow H\ot C,\;\;\;\lambda (c)=c^{(-1)}\ot c^{(0)},
\end{eqnarray*}
for all $c\in C$. We record the (bi) comodule conditions:
\begin{eqnarray}
&&c^{(-1)}\ot c^{(0)(-1)}\ot c^{(0)(0)}=(c^{(-1)})_1\ot (c^{(-1)})_2
\ot c^{(0)}, \label{lc} \\
&&c^{<0><0>}\ot c^{<0><1>}\ot c^{<1>}=c^{<0>}\ot (c^{<1>})_1\ot 
(c^{<1>})_2, \label{rc} \\
&&c^{<0>(-1)}\ot c^{<0>(0)}\ot c^{<1>}=c^{(-1)}\ot c^{(0)<0>}\ot c^{(0)<1>}.
\label{bc}
\end{eqnarray}
(ii) $C$ is a coalgebra, with comultiplication $\Delta _C:C\rightarrow 
C\ot C$, $\Delta _C(c)=c_1\ot c_2$, and counit $\varepsilon _C:C\rightarrow 
k$.\\
(iii) $C$ is a left $H$-comodule coalgebra, that is, for all $c\in C$:
\begin{eqnarray}
&&c_1^{(-1)}c_2^{(-1)}\ot c_1^{(0)}\ot c_2^{(0)}=c^{(-1)}\ot 
(c^{(0)})_1\ot (c^{(0)})_2, \label{lca} \\
&&c^{(-1)}\varepsilon _C(c^{(0)})=\varepsilon _C(c)1_H.
\end{eqnarray}
(iv) $C$ is a right $H$-comodule coalgebra, that is, for all $c\in C$:
\begin{eqnarray}
&&c_1^{<0>}\ot c_2^{<0>}\ot c_1^{<1>}c_2^{<1>}=(c^{<0>})_1\ot (c^{<0>})_2
\ot c^{<1>}, \label{rca} \\
&&\varepsilon _C(c^{<0>})c^{<1>}=\varepsilon _C(c)1_H.
\end{eqnarray}
We denote $C\ot H$ by $C\nat H$ and elements $c\ot h$ by $c\nat h$.  
Define the maps 
\begin{eqnarray*}
&&\Delta :C\nat H\rightarrow (C\nat H)\ot (C\nat H), \;\;\varepsilon :
C\nat H\rightarrow k, \\ 
&&\Delta (c\nat h)=(c_1^{<0>}\nat c_2^{(-1)}h_1)\ot (c_2^{(0)}\nat 
h_2c_1^{<1>}), \;\; 
\varepsilon (c\nat h)=\varepsilon _C(c)\varepsilon _H(h). 
\end{eqnarray*}
\begin{proposition}
$(C\nat H, \Delta , \varepsilon )$ is a coalgebra, called the  
L-R-smash coproduct. 
\end{proposition}
\begin{proof} 
The counit axiom is immediate, so we check coassociativity. We compute:
\begin{eqnarray*}
(id\ot \Delta )(\Delta (c\nat h))&=&
(c_1^{<0>}\nat c_2^{(-1)}h_1)\ot \Delta (c_2^{(0)}\nat h_2c_1^{<1>})\\
&=&(c_1^{<0>}\nat c_2^{(-1)}h_1)\ot ((c_2^{(0)})_1^{<0>}\nat 
(c_2^{(0)})_2^{(-1)}h_2(c_1^{<1>})_1)\\
&&\ot ((c_2^{(0)})_2^{(0)}\nat h_3(c_1^{<1>})_2(c_2^{(0)})_1^{<1>})\\
{\rm (\ref{lca})}&=&(c_1^{<0>}\nat c_2^{(-1)}c_3^{(-1)}h_1)\ot 
(c_2^{(0)<0>}\nat c_3^{(0)(-1)}h_2(c_1^{<1>})_1)\\
&&\ot (c_3^{(0)(0)}\nat h_3(c_1^{<1>})_2c_2^{(0)<1>})\\
{\rm (\ref{rc})}&=&(c_1^{<0><0>}\nat c_2^{(-1)}c_3^{(-1)}h_1)\ot 
(c_2^{(0)<0>}\nat c_3^{(0)(-1)}h_2c_1^{<0><1>})\\
&&\ot (c_3^{(0)(0)}\nat h_3c_1^{<1>}c_2^{(0)<1>})\\
{\rm (\ref{bc})}&=&(c_1^{<0><0>}\nat c_2^{<0>(-1)}c_3^{(-1)}h_1)\ot 
(c_2^{<0>(0)}\nat c_3^{(0)(-1)}h_2c_1^{<0><1>})\\
&&\ot (c_3^{(0)(0)}\nat h_3c_1^{<1>}c_2^{<1>})\\
{\rm (\ref{lc})}&=&(c_1^{<0><0>}\nat c_2^{<0>(-1)}(c_3^{(-1)})_1h_1)\ot 
(c_2^{<0>(0)}\nat (c_3^{(-1)})_2h_2c_1^{<0><1>})\\
&&\ot (c_3^{(0)}\nat h_3c_1^{<1>}c_2^{<1>})\\
{\rm (\ref{rca})}&=&((c_1^{<0>})_1^{<0>}\nat (c_1^{<0>})_2^{(-1)}
(c_2^{(-1)})_1h_1)\ot ((c_1^{<0>})_2^{(0)}\nat 
(c_2^{(-1)})_2h_2(c_1^{<0>})_1^{<1>})\\
&&\ot (c_2^{(0)}\nat h_3c_1^{<1>})\\
&=&\Delta (c_1^{<0>}\nat c_2^{(-1)}h_1)\ot (c_2^{(0)}\nat h_2c_1^{<1>})\\
&=&(\Delta \ot id)(\Delta (c\nat h)), 
\end{eqnarray*}
finishing the proof.
\end{proof}
\begin{remark}
Obviously, if the right $H$-comodule structure of $C$ is trivial (i.e. 
$c^{<0>}\ot c^{<1>}=c\ot 1$ for all $c\in C$), then $C\nat H$ coincides 
with Molnar's smash coproduct from \cite{molnar}.
\end{remark}

\end{document}